\documentclass{article}

\usepackage{arxiv}

\usepackage[utf8]{inputenc} 
\usepackage[T1]{fontenc}    
\usepackage{url}            
\usepackage{booktabs}       
\usepackage{amsfonts}       
\usepackage{microtype}      
\usepackage{lipsum}
\usepackage{amssymb}
\usepackage{amsmath}
\usepackage{amscd}
\usepackage{graphicx}
\usepackage{graphics}
\usepackage[noadjust]{cite}
\usepackage{amsthm}
\usepackage{caption}
\usepackage{multirow}
\usepackage{array}
\usepackage{rotating}
\usepackage[norelsize,ruled,linesnumbered]{algorithm2e}
\usepackage{algorithm2e,setspace}
\usepackage{hyperref}

\usepackage{indentfirst}
\usepackage{tikz}
\usepackage{calc}
\usepackage{epsfig}
\usepackage{natbib}
\usepackage[makeroom]{cancel}

\title{Improving full  waveform inversion by wavefield reconstruction with the alternating direction 
method of multipliers}

\author{
  Hossein S.~Aghamiry\\
  Institute of Geophysics, University of Tehran, Tehran, Iran.
  \texttt{h.aghamiry@ut.ac.ir} \\
  Geoazur, Universit\'e C\^ote d'Azur, CNRS, IRD, OCA, Valbonne, France. 
  \texttt{aghamiry@geoazur.unice.fr} 
   \And
 Ali Gholami \\
  Institute of Geophysics, University of Tehran, Tehran, Iran.
  \texttt{agholami@ut.ac.ir} 
  \And
  St\'ephane Operto \\
  Geoazur, Universit\'e C\^ote d'Azur, CNRS, IRD, OCA, Valbonne, France. 
  \texttt{operto@geoazur.unice.fr}
  }

\begin{document}
\maketitle

\begin{abstract}
Full waveform inversion (FWI) is an iterative nonlinear waveform matching procedure subject to wave-equation constraint. FWI is highly nonlinear when the wave-equation constraint is enforced at each iteration. To mitigate  nonlinearity, wavefield-reconstruction inversion (WRI) expands the search space by relaxing the wave-equation constraint with a penalty method. The pitfall of this approach resides in the tuning of the penalty parameter because increasing values should be used to foster data fitting during early iterations while progressively enforcing the wave-equation constraint during late iterations. However, large values of penalty parameter lead to ill-conditioned problems. Here, this tuning issue is solved by replacing the penalty method by an augmented Lagrangian method equipped with operator splitting (IR-WRI as iteratively-refined WRI).  It is shown that IR-WRI is similar to a penalty method in which data and sources are updated at each iteration by the running sum of the data and source residuals of previous iterations. Moreover, the alternating direction strategy exploits the bilinearity of the wave equation constraint to linearize the subsurface model estimation around the reconstructed wavefield. Accordingly, the original nonlinear FWI is decomposed into a sequence of two linear subproblems, the optimization variable of one subproblem being passed as a passive variable for the next subproblem. The convergence of WRI and IR-WRI are first compared with a simple transmission experiment, which lies in the linear regime of FWI. Under the same conditions, IR-WRI converges to a more accurate minimizer with a smaller number of iterations than WRI. More realistic case studies performed with the Marmousi II and the BP salt models show the resilience of IR-WRI to cycle skipping and noise, as well as its ability to reconstruct with high fidelity large-contrast salt bodies and sub-salt structures starting the inversion from crude initial models and a 3-Hz starting frequency.
\end{abstract}


\section{Introduction}
Full Waveform Inversion (FWI) is a leading-edge seismic imaging method that has gained renewed interest in oil and gas exploration since one decade due to the emergence of supercomputers, long-offset wide-azimuth acquisition, and broadband seismic sources \citep{Virieux_2009_OFW,Virieux_2017_FWI}.
FWI seeks to estimate subsurface parameters 
with a resolution close to the seismic wavelength by fitting seismic records with modeled seismograms that are computed with two-way modeling engines to account for the full complexity of the wavefields \citep{Tarantola_1984_ISR,Pratt_1998_GNF}. 
%

From the numerical optimization viewpoint, FWI can be defined as a partial differential equation (PDE)-constrained nonlinear optimization problem, where the equality constraint is the wave equation.
A popular method to solve such nonlinear constrained optimization problem is the method of Lagrange multiplier, where the constraints are implemented through a Lagrangian function \citep[e.g.][]{Haber_2000_OTS}. A local minimizer satisfies the so-called KKT (Karush-Kuhn-Tucker) first-order optimality conditions at the stationary point of the Lagrangian. 
In the frame of local optimization methods, the so-called full-space approach satisfies the nonlinear KKT conditions  with an iterative Newton method which jointly updates three classes of variables (the state variables, the adjoint-state variables and the parameters). However, the dimension of the KKT system makes this full space approach computationally intractable for FWI application.
Instead, FWI is classically solved with a reduced approach, namely a variational projection approach where the full-space optimization problem is recast as an unconstrained problem after successive elimination of the state and adjoint-state variables from the KKT system \citep{Haber_2000_OTS,Akcelik_2002_MNK,Askan_2007_FWI,Epanomeritakis_2008_NCG}. 
Although this reduced-approach strictly enforces the constraints (because the wave equation is solved exactly with the current guess of the subsurface parameters) and is computationally tractable, it is highly nonlinear and ill-posed due to the oscillating nature of seismic signals \citep[e.g. ][]{Symes_2008_MVA}.
One of the main source of nonlinearity is cycle skipping, which happens as soon as the current subsurface model is not accurate enough to predict recorded traveltimes with an error lower than half the dominant period \citep[e.g.][ Figure 7]{Virieux_2009_OFW}. In contrast, full-space approaches are more resilient to cycle skipping because the joint updating of wavefields, adjoint wavefields and subsurface model makes them more versatile to fit the data with inaccurate subsurface models.
%
%
The cycle skipping criterion is difficult to satisfy when the number of propagated wavelengths increases as might be expected with the development of long-offset wide-azimuth acquisitions \citep{Pratt_2008_WAT,Virieux_2009_OFW}. This prompts the oil industry to develop broadband sources with a richer low frequency content to balance increasing propagation distances and mitigate cycle skipping accordingly \citep{Plessix_2012_FWI,Baeten_2013_ULF}. In parallel with this, many heuristic data-driven approaches based upon frequency, traveltime and offset continuations have been designed to drive the optimization problem toward multiscale imaging from the long wavelengths to the shorter ones \citep[e.g.][]{Gorszczyk_2017_FWI}. To reduce human intervention, more automatic approaches rely on more convex distances between recorded and modeled data such as those based upon instantaneous phase and envelope \citep{Bozdag_2011_MFF,Luo_2015_ENV} or correlation \citep{VanLeeuwen_2010_CMC}. However, these approaches should be rather viewed as tools to build initial models for FWI rather than FWI methods per se because they often achieve an improved convexity at the expense of spatial resolution.

More ambitious approaches extend the linear regime of FWI by enlarging the search space. 
Extended-domain approaches have been intensively developed for velocity macromodel building in the framework of migration-based velocity analysis (\citet{Symes_2008_MVA} for a review). These approaches extends the model space with additional degree of freedom as subsurface offset or time lag, which makes them quite computationally intensive. Then, the inverse problem seeks a velocity model that produces a physical reflectivity through an annihilator applied to the extended migration operator. More recently, \citet{VanLeeuwen_2013_MLM} extended the search space in frequency-domain FWI by replacing the wave-equation constraint by a quadratic penalty term. This penalty method, called wavefield-reconstruction inversion (WRI), relaxes the requirement to satisfy exactly the wave equation at each iteration for the benefit of an improved data fitting. This is achieved by reconstructing the wavefield that best jointly fits in the least-squares sense the observations and satisfy the wave equation.
To make WRI computationally tractable, \citet{VanLeeuwen_2013_MLM} perform the wavefield reconstruction and the subsurface parameter estimation in an alternating way: the wavefield is first estimated using the available subsurface model as fixed background model, before updating the subsurface parameters using the previously-reconstructed wavefield as a fixed background wavefield, this cycle being iterated until convergence. This alternating-direction strategy linearizes the parameter-estimation subproblem around the reconstructed wavefield. Accordingly, parameter updating can be performed with Gauss-Newton iterations, the Hessian being diagonal when the forward problem equation is the Helmholtz equation \citep[][ their equation 8]{VanLeeuwen_2013_MLM}. Later, \citet{vanLeeuwen_2016_PMP} reformulated WRI as a reduced penalty method implemented with a variable projection approach: the closed-form expression of the extended-domain reconstructed wavefield is injected as a function of the subsurface parameters in the penalty function instead of using this wavefield as a passive variable (i.e., independent to the subsurface parameters). Although this variable elimination makes the parameter-estimation subproblem non linear, \citet{vanLeeuwen_2016_PMP} assess their method with a Gauss-Newton method (by opposition to the full Newton counterpart) to mitigate the computational burden. Moreover, using a sparse approximation of the Gauss-Newton Hessian makes the descent direction of the reduced approach identical to that of the alternating-direction WRI of \citet{VanLeeuwen_2013_MLM}. In the following, for sake of concise notations, we will refer WRI to as the alternating-direction penalty method of \citet{VanLeeuwen_2013_MLM}, because it is more closely related to the approach developed in this study.
Compared to full-space approaches, WRI is less computationally expensive as each source is processed independently and the parameter estimation does not require solving an adjoint-state equation because the Lagrange multiplier has been eliminated from the optimization problem at the benefit of a penalty parameter. \\
%
One potential difficulty with WRI resides in the tuning of the penalty parameter $\lambda$. Ideally, increasing values should be used during iterations to  progressively enforce the wave-equation constraint and, hence satisfy the KKT optimality conditions of the original constrained problem with acceptable precision at the minimizer. A significant issue with this continuation approach is that the Hessian is ill conditioned for large $\lambda$. However, \citet{vanLeeuwen_2016_PMP} suggested that a fixed $\lambda$ equal to a fraction of the largest eigenvalue of the normal operator associated with the augmented wave-equation system may allow for convergence toward solutions with acceptable error bounds. \\
%
%
For the sake of completeness, several other extended-search methods, which share the idea of fitting the data by violating the wave equation constraint and updating the subsurface model by minimizing the constraint violation, have been proposed.  \citet{Luo_2011_DBO,Warner_2016_AWI,Zhu_2016_BSM} estimate the filter that allows for the matching between the recorded and modeled data before minimizing the deviation of this matching filter from a delta function. Closely related to this approach, the source-receiver extension method first matches the data through a trace-per-trace source signature estimation before updating the subsurface model by minimizing the deviation of each non-physical source signature from a reference one through an anihilator \citep{Huang_2017_FWI}. In line with the approach promoted by \citet{Symes_2008_MVA}, \citet{Biondi_2014_SIF} extend the velocity model along the time-lag axis to account for large traveltime residuals and augment the misfit function with a penalty term which penalizes large time lags through an anihilator. With the same objective of managing large time shifts, \citet{Engquist_2014_WAS,Metivier_2015_MTM,Metivier_2016_OTI,Yang_2017_AOT} compute the optimal transport between the observation and the modeled data before updating the subsurface parameters by minimizing the optimal-transport distance. \citet{Ma_2013_WRT} build traveltime residual maps between recorded and modeled shot gathers by dynamic warping and minimize in a second step these traveltime residuals. 

%
%
In this framework, this study focuses on an easy-to-implement improvement of WRI, which strongly reduces the sensitivity of the method to the penalty parameter, resulting in a faster convergence and more accurate subsurface and wavefield reconstruction.
%
%
To achieve this goal, we replace the penalty method by an augmented Lagrangian method which is implemented with a classical primal descent - dual ascent approach. In constrained optimization, the augmented Lagrangian method can be viewed as a combination of a classical Lagrangian method and a penalty method, where the Lagrangian multiplier steers the inversion toward the constraint, while the penalization term regularizes the inversion \citep[][ Chapter 17]{Nocedal_2006_NOO}. As in WRI, we perform wavefield reconstruction and parameter estimation in an alternating mode to make the computational problem tractable, while the Lagrange multipliers are updated with a classical dual ascent approach.
This operator splitting draw some connection between our approach and the popular alternating direction method of multiplier (ADMM) \citep{Boyd_2004_COO} with the difference that the wavefield reconstruction and the parameter estimation subproblems are not separable. This non-separability is managed by solving the two subproblems in sequence rather than in parallel, the solution of one subproblem being passed as a passive variable to the next problem. As above mentioned, this sequential splitting strategy linearizes the parameter estimation problem around the reconstructed wavefield, which means that the original non linear FWI has been recast as a sequence of two linear subproblems. \\
By using a scaled form of the Lagrangian, we also show that the augmented Lagrangian method is similar to the WRI penalty method, except that the right-hand sides in the two objectives of the penalty function (namely, the data and the seismic sources) are updated at each iteration with the running sum of the data and source residuals of previous iterations. This right-hand side updating decreases more efficiently the data and source residuals in iterations and performs a self-adaptive weighting of the competing data-fitting and wave-equation objectives of the penalty function, leading to more accurate minimizers and faster convergence. Accordingly, we refer our method to as IR-WRI (Iteratively-Refined WRI). \\
We assess our method against three synthetic examples. First, a toy example highlights how the convergence history of IR-WRI relative to WRI is improved by data and source residual updating.
Then, we apply IR-WRI to two challenging synthetic models (the Marmousi II and the 2004 BP salt models) starting from crude initial models and realistic frequency (3~Hz). As WRI, IR-WRI shows a good resilience 
to cycle skipping. However, IR-WRI converges to subsurface models of much improved resolution and, in addition, it is more resilient to noise.

This paper is mainly organized in a method, numerical example and discussion sections. In the method section, we first review the principles of FWI based on full-space and reduced-space approaches before introducing WRI. Then, we describe IR-WRI. We first introduce a scaled form of the augmented Lagrangian method to highlight the similarities and differences between WRI and IR-WRI. Finally, we review the implementation of the operator splitting strategy, in particular the dual updating of the Lagrange multipliers. In the second part of the paper, we present three numerical examples where we compare the results of WRI and IR-WRI.
We first illustrate the impact of the right-hand side updating (data and source) on the convergence history of the method with a toy example. Then, we illustrate the ability of IR-WRI to image large-contrast media such as the BP salt model. In the final discussion section, the results of the numerical examples are interpreted in light of the IR-WRI methodological ingredients, some analogies with other optimization method such as the Bregman method are discussed and some perspective works are introduced.
%

%
%
\section{Method}

\subsection*{FWI with full-space and reduced-space Lagrange multiplier method}

The FWI can be formulated in the frequency domain as the following nonlinear multivariate PDE-constrained optimization problem \citep[e.g.][]{VanLeeuwen_2013_MLM,vanLeeuwen_2016_PMP}:
\begin{equation} \label{orginal_FWI}
\min_{\bold{m},\bold{u}} ~~~~~~ \|\bold{Pu}-\bold{d}\|_2^2, ~~~~ \text{subject to} ~~~~ \bold{A(m)u}=\bold{b},
\end{equation}
where $\|\cdot\|_2^2$ denotes the Euclidean norm, $\bold{m} \in \mathbb{R}^{N\times 1}$ the discrete subsurface model parameters,
$\bold{b} \in \mathbb{C}^{N\times 1}$ the source term,
$\bold{u} \in \mathbb{C}^{N\times 1}$ the modeled wavefield,
$\bold{d} \in \mathbb{C}^{M\times 1}$ the recorded seismic data and $\bold{P} \in \mathbb{R}^{M\times N}$ a linear observation operator that samples $\bold{u}$ at the receiver positions. 
We consider one frequency and one source in equation~\ref{orginal_FWI} for sake of notation compactness. Extension to multiple frequencies and sources is simply implemented by 
summation over sources and frequencies in the objective and by adding multiple right-hand sides in the constraint, equation ~\ref{orginal_FWI}.
The matrix $\bold{A(m)} \in \mathbb{C}^{N\times N}$, whose coefficients depend on $\bold{m}$, represents the discretized PDE \citep{Pratt_1998_GNF,Plessix_2007_HIS}.
In this study, we limit ourselves to the Helmholtz equation assuming a constant density equal to 1.
\begin{eqnarray}  
\bold{A(m)} &= & \bold{\Delta} + \omega^2 \bold{B}\bold{C}(\bold{m}) \text{diag}(\bold{m}),
\label{helmholtz}
\end{eqnarray}
where $\omega$ is the angular frequency, $\bold{\Delta}$ a discretized Laplace operator and the subsurface parameters $\bold{m}$ are parametrized by the squared slowness. 
The operator $\bold{C}$ enclose boundary conditions, which can be a function of $\bold{m}$ (e.g., Robin paraxial conditions \citep{Engquist_1977_ABC}) or independent from $\bold{m} $ (e.g., sponge-like absorbing boundary conditions such as perfectly-matched layers \citep{Berenger_1994_PML}). Also, the linear operator $\bold{B}$ spreads the "mass" term $\omega^2\bold{C}(\bold{m}) \text{diag}(\bold{m})$ 
over all the  coefficients of the stencil to improve its accuracy following an anti-lumped mass strategy \citep{Marfurt_1984_AFF,Jo_1996_OPF,Hustedt_2004_MGS}.   

The nonlinear constrained optimization problem, equation \ref{orginal_FWI}, can be solved with the method of Lagrange multipliers and Newton-type algorithm \citep{Haber_2000_OTS}. 
\begin{equation} \label{orginal_FWI1}
\min_{\bold{m},\bold{u}} \max_{\bold{v}} \mathcal{L} (\bold{m},\bold{u},\bold{v}) = \min_{\bold{m},\bold{u}} \max_{\bold{v}} \|\bold{Pu}-\bold{d}\|_2^2 + \bold{v}^T \left[ \bold{A(m)u} - \bold{b} \right],
\end{equation}
where $\bold{v} \in \mathbb{C}^{N\times 1}$ denotes the Lagrange multiplier (or, adjoint-state wavefield in the FWI terminology) and $^T$ the complex-conjugate (Hermitian) transpose.

A full-space approach jointly updates the three classes of variables ($\bold{u}$, $\bold{v}$, $\bold{m}$). A Newton approach would require to solve the large-scale KKT system, whose matrix operator is the multivariate Hessian, the unknowns are gathered in the multivariate descent direction and the right-hand side is the multivariate gradient of the Lagrangian \citep{Haber_2000_OTS,Akcelik_2002_MNK,Askan_2007_FWI}.
Instead, a reduced approach, which strictly enforces the PDE constraint at each iteration by projection of the full multivariate search space onto the parameter search space, is more commonly use for sake of computational efficiency. This projection can be implemented by successive elimination of the state and adjoint-state variables from the KKT system \citep{Akcelik_2002_MNK} or, more simply, by enforcing the wave-equation solution in 
the objective leading to the following unconstrained optimization problem \citep{Pratt_1998_GNF,Plessix_2006_RAS}
\begin{equation} \label{pratt_FWI}
\min_{\bold{m}} J_r(\bold{m}) = \min_{\bold{m}}  \|\bold{P}\bold{A}^{-1}(\bold{m})\bold{b}-\bold{d}\|_2^2.
\end{equation}
The gradient of the reduced objective function at iteration $k+1$ is 
the zero-lag correlation between the so-called virtual scattering source, namely the reconstructed incident wavefield, $\bold{u}^{k+1}=\bold{A}^{-1}(\bold{m}^k) \bold{b}$, weighted by the radiation pattern $\partial \bold{A}(\bold{m}^k) / \partial \bold{m}$ \citep{Pratt_1998_GNF} and $\bold{v}^{k+1}$ 
\begin{equation}
\nabla_{\bold{m}} J_r = \Re \left(  \omega^2 \text{diag}(\bold{u}^{k+1} ) ^T \bold{v}^{k+1} \right),
\end{equation}
where the radiation pattern of the squared slowness simply samples $\bold{u}^{k+1}$ at the position of the model parameter with a weight $\omega^2$ and
$\bold{v}^{k+1}$ satisfies the so-called adjoint-state equation
\begin{equation}
\bold{A}^T(\bold{m}^k) \bold{v}^{k+1} = \bold{P}^T \big[ \bold{P} \bold{u}^{k+1} - \bold{d} \big],
\label{eqadjoint}
\end{equation}
whose right-hand side is the prolongation of the data residuals in the full computational domain. \\
Although this reduced approach requires to solve exactly the state and adjoint-state equations at each nonlinear iteration, it is more computationally tractable than the full-space approach. Moreover, after elimination of the state and adjoint-state variables, one needs to update only one class of unknowns, namely $\bold{m}$, at each iteration.
However, due to the highly-oscillating nature of  the inverse PDE operator $\bold{A}^{-1}$ (the Green functions), the elimination of variables associated with the reduced approach makes the inverse problem highly nonlinear, and hence prone to convergence to inaccurate minimizer when the initial $\bold{m}$ is not accurate enough \citep{Symes_2008_MVA,Virieux_2009_OFW}.

\subsection*{Wavefield-Reconstruction Inversion (WRI): alternating-direction penalty method}

To overcome the limitations of the reduced approach, \citet{VanLeeuwen_2013_MLM} recast the original constrained optimization problem, equation \ref{orginal_FWI}, as a multi-variate unconstrained quadratic penalty problem for $\bold{u}$ and $\bold{m}$ given by
\begin{equation} \label{leeuwen_FWI}
\min_{\bold{m,u}} J_p(\bold{m,u}) =  \min_{\bold{m,u}} \|\bold{Pu}-\bold{d}\|_2^2 + \lambda \|\bold{A(m)u}-\bold{b}\|_2^2,
\end{equation}
where the scalar $\lambda > 0$ is the so-called penalty parameter. The hard constraint in the original constrained problem is replaced by a quadratic penalty term, which represents the squares of the constraint violation.
Also, the penalty problem, equation~\ref{leeuwen_FWI}, can be viewed as a compromise between the full-space and the reduced approach in the sense that the full search space $(\bold{m},\bold{u},\bold{v})$ has been projected 
onto the $(\bold{m},\bold{u})$ space by enforcing $\bold{v} = \lambda (\bold{A(m)u}-\bold{b})$ in the Lagrangian function, equation~\ref{orginal_FWI1} \citep{vanLeeuwen_2016_PMP}.
By allowing for the PDE constraint violation, \citet{VanLeeuwen_2013_MLM} enlarge the search space and potentially mitigate the inversion nonlinearity.
However, the joint update of $\bold{u}$ and $\bold{m}$ remains challenging from the computational viewpoint and \citet{VanLeeuwen_2013_MLM} resort to an alternating-direction 
minimization strategy to solve the penalty problem.
A cycle of the algorithm first reconstructs $\bold{u}$ by solving in a least-squares sense the overdetermined system for each source $\bold{b}$ 
\begin{equation}
\begin{bmatrix}
\sqrt{\lambda} {\bf{A}} \\
\bold{P} \\
\end{bmatrix}
\bold{u}
=
\begin{bmatrix}
\sqrt{\lambda} \bold{b} \\
\bold{d}
\end{bmatrix}.
\label{system}
\end{equation}
Equation \ref{system} relaxes the PDE constraint, while pushing the wavefield reconstruction against the observations through the feedback term $(\bold{P} \bold{u} = \bold{d})$.  
At iteration $k+1$ of the alternating-direction workflow, the closed-form expression of $\bold{u}$ \citep[][ equation 19]{vanLeeuwen_2016_PMP} is given by
\begin{equation}
\bold{u}^{k+1}=\left[\bold{A}^T(\bold{m}^{k}) \bold{A}(\bold{m}^{k})  + \lambda^{-1} \bold{P}^T \bold{P}\right]^{-1} \left[\bold{A}^T(\bold{m}^{k})  \bold{b} + \lambda^{-1} \bold{P} \bold{d} \right],
\label{equvlh}
\end{equation}
where the update of $\bold{m}$ at the previous iteration, namely $\bold{m}^{k}$, is used as a "passive" variable according to the 
alternating-direction strategy. Then, $\bold{u}^{k+1}$ is passed as a passive variable to the second subproblem for $\bold{m}$. Assuming $\bold{u}$ fixed (namely, independent to $\bold{m}$) makes the subproblem for $\bold{m}$ linear \citep[][ their equation 8]{VanLeeuwen_2013_MLM}. Accordingly, $\bold{m}$ is updated with one Gauss-Newton iteration before proceeding with the next cycle of the alternating-direction optimization.
At iteration $k+1$, the gradient of $J_p$ with respect to $\bold{m}$ is simply formed by the zero-lag correlation between the virtual scattering source $\left(\partial \bold{A}(\bold{m})^{k} / \partial \bold{m}\right)  \bold{u}^{k+1}$ and the source residual $\bold{A(m}^{k})\bold{u}^{k+1}-\bold{b}$, \citep[][ equation 7]{VanLeeuwen_2013_MLM}, namely
\begin{equation}
\nabla_\bold{m} J_p=\Re \left(   \omega^2 \text{diag}(\bold{u}^{k+1})^T \big[ \bold{A}(\bold{m}^{k}) \bold{u}^{k+1} - \bold{b}\big] \right).
\label{gradvlgh}
\end{equation}
\citet{vanLeeuwen_2016_PMP} formulated a variant of the alternating-direction penalty method based upon a reduced variable projection approach. The closed-form expression of the extended-space $\bold{u}$ as a function of $\bold{m}$, equation~\ref{equvlh} where the fixed $\bold{m}^{k}$ should be replaced by the variable $\bold{m}$, is injected in the penalty function before updating $\bold{m}$ with a nonlinear Newton iteration. The gradient of the reduced penalty function at $\bold{m}=\bold{m}^k$ has the same expression as in equation~\ref{gradvlgh}, but the Hessian embeds more complex first-order and second-order terms resulting from the dependency of $\bold{u}$ on $\bold{m}$, \citep[][ Their equations 14-17]{vanLeeuwen_2016_PMP}. \citet{vanLeeuwen_2016_PMP} neglect second-order terms and update $\bold{m}$ with Gauss-Newton iterations for sake of computational efficiency. \citet{vanLeeuwen_2016_PMP} also show that a sparse Gauss-Newton approximation of the reduced approach is equivalent to the original alternating-direction WRI of \citet{VanLeeuwen_2013_MLM}. From now on, we will refer WRI to as the alternating-direction implementation of the penalty method to which IR-WRI is more closely related.

As mentioned in the introduction, a well known limitation of penalty methods is the tedious tuning of the penalty parameter through iterations. Ideally, a sequence of increasing $\lambda_k$ should be used such that the constraint is satisfied with a prescribed tolerance at the last iteration of the minimization problem. However, using large values of $\lambda_k$ generally leads to ill-conditioned and unstable problems \citep[][ page 505]{Nocedal_2006_NOO}. 
Moreover, one may have noted that, at a given iteration of the wavefield reconstruction, equation ~\ref{equvlh}, the right-hand sides $\bold{d}$ and  $\bold{b}$ of the two objectives, equation~\ref{leeuwen_FWI}, are not updated by the running sum of the errors (or, residuals) of previous iterations, although updating of right-hand sides is a common practice to iteratively refine the solution of linear inverse problems (Appendix 1). The same comment applies to the parameter estimation subproblem in \citep[][ their equation 8]{VanLeeuwen_2013_MLM}, where the source is not updated with the residuals of previous iterations. This might prevent an efficient decrease of the data residuals and  wave-equation error in iterations.

In the next section, we review our approach based on the augmented Lagrangian method with operator splitting to overcome these two limitations.



\subsection*{WRI with an alternating-direction method of multipliers}
\subsubsection*{Feasibility problem and scaled-form augmented Lagrangian}
Let us recast FWI as the following feasibility problem
\begin{equation}
\label{main}
{\text{find}} ~~~~ \bold{m} ~\text{and}~ \bold{u} ~~~~\text{subject to} ~~~~\bold{A(m)u}=\bold{b} ~~~~  \text{and}    ~~~~   \bold{Pu}= \bold{d},
\end{equation} 
where the term $\bold{Pu} = \bold{d}$ is now processed as an equality constraint in case of noiseless data. In presence of noise, this equality constraint becomes an inequality constraint, $\|\bold{Pu} - \bold{d}\|_2^2 < \epsilon_n$, where $\epsilon_n$ represents the noise level. We assume that the feasible set is nonempty, namely the constraints are consistent or there are vectors $\bold{m}$ and $\bold{u}$ satisfying the constraints. The feasibility problem can be considered as a special case of a constrained optimization problem where the objective function is identically zero \citep[][ pages 128-129]{Boyd_2004_COO}. \\
Accordingly, the feasibility problem, equation~\ref{main}, is recast as the constrained problem
\begin{equation}
\label{main1_1}
\min_{\bold{u},\bold{m}} ~~ \bold{0} ~~~~ \text{subject to} ~~~~ \bold{A(m)u}=\bold{b} ~~~~ \text{and} ~~~~  \bold{Pu}= \bold{d}.  
\end{equation} 
where the first and second constraints are bilinear and linear, respectively. 
The non-linear feasibility problem \ref{main1_1} can be written in more compact form as
\begin{equation}
\label{main1_12}
\min_{\bold{u},\bold{m}} ~~ \bold{0} ~~~~ \text{subject to} ~~~~ \bold{F(m)u}=\bold{s}, 
\end{equation} 
where
\begin{equation}
 \bold{s} =
 \begin{bmatrix}
\bold{d} \\
\bold{b}
\end{bmatrix} \in \mathbb{C}^{(M+N)\times 1}
~~~~~ \text{and} ~~~~~
 \bold{F}(\bold{m}) =
 \begin{bmatrix}
\bold{P} \\
\bold{A(m)}
\end{bmatrix} \in \mathbb{C}^{(M+N)\times N}. \nonumber
\end{equation}
The augmented Lagrangian function, \citet[][ their equation 2.6]{Boyd_2011_DOS} and \citet[][ their equation 17.36]{Nocedal_2006_NOO}, for the problem defined by equation \ref{main1_12} is 
\begin{equation} \label{AL}
\mathcal{L}(\bold{m},\bold{u},\bold{v}) = \bold{v}^T[\bold{F(m)u}-\bold{s}]  + \frac{1}{2}\|\bold{\Gamma}^\frac12\bold{[F(m)u}-\bold{s}]\|_2^2,
\end{equation}
where
\begin{equation} 
\bold{\Gamma} =
 \left(
    \begin{array}{r@{}c|c@{}l}
  &    \begin{smallmatrix}
       \gamma_0 & & 0 \\
          &\ddots&\\
        0 & & \gamma_0 \rule[-1ex]{0pt}{2ex}
      \end{smallmatrix} 
  & \mbox{\huge0} 
  & \rlap{\kern5mm$M$}\\
      \hline
  &    \mbox{\huge0} &  
       \begin{smallmatrix}\rule{0pt}{2ex}
        \gamma_1 & & 0 \\
          &\ddots&\\
        0 & & \gamma_1
      \end{smallmatrix}    &  \rlap{\kern5mm $N$}
    \end{array} 
\right), \nonumber
\end{equation}
%
$\gamma_0$, $\gamma_1 >0$ are penalty parameters and $\bold{v} \in \mathbb{C}^{(M+N)\times 1}$ is the Lagrangian multiplier (known as dual variable).
Applying the method of multipliers, \citet[][ their equations 2.7-2.8]{Boyd_2011_DOS} and \citet[][ their equation 17.39]{Nocedal_2006_NOO}, we get
\begin{subequations}
\label{main1}
\begin{eqnarray}
&&(\bold{m}^{k+1},\bold{u}^{k+1}) = \underset{\bold{m},\bold{u}}{\arg\min} ~~ \mathcal{L}(\bold{m},\bold{u},\bold{v}^k), \label{main1a} \\
&&\bold{v}^{k+1} = \bold{v}^{k} +  \bold{\Gamma}[\bold{F}(\bold{m}^{k+1})\bold{u}^{k+1}-\bold{s}],  \label{main1c}
\end{eqnarray} 
\end{subequations}
for $k=0,1,...$ beginning from $\bold{v}^0=\bold{0}$. 
This iteration can be interpreted as follows:
we begin with a prior estimate $\bold{v}^0$ of the dual variable $\bold{v}$ and minimize the objective function with respect to the primal variables $\bold{m}$ and $\bold{u}$, equation \ref{main1a}. Then, we update the dual variable $\bold{v}$ via a steepest ascent method by keeping the primal variables $\bold{m}$ and $\bold{u}$ fixed (equation \ref{main1c}).
This process is iterated until convergence, i.e., when $\bold{F}(\bold{m}^{k+1})\bold{u}^{k+1}=\bold{s}$.

To highlight more explicitly the key role of iterative right-hand side updating in our method and the relationship with the penalty method of \citet{VanLeeuwen_2013_MLM}, we rewrite the Lagrangian, equation~\ref{AL}, in a scaled form \citep[][ Section 3.1.1]{Boyd_2011_DOS}. 
Introducing the scaled dual variable $\bold{s}^k=-\bold{\Gamma}^{-1}\bold{v}^k$ in the Lagrangian and adding the constant terms $\frac{1}{2}\|\Gamma^\frac12\bold{s}^k\|_2^2$ to it, we find
\begin{equation} \label{ALS0}
\mathcal{L}(\bold{m},\bold{u},\bold{s}^k) = \frac{1}{2}  \|\Gamma^\frac12[\bold{F(m)u}-\bold{s}-\bold{s}^k]\|_2^2.
\end{equation}
Without modifying its minimum, we can scale the augmented Lagrangian \ref{ALS0}  by $1/\gamma_0$ to eliminate one of the penalty parameters and get
\begin{equation} \label{ALS}
\mathcal{L}(\bold{m},\bold{u},\bold{s}^k) = \frac{1}{2}  \|\bold{\Lambda}^\frac12[\bold{F(m)u}-\bold{s}-\bold{s}^k]\|_2^2,
\end{equation}
where 
\begin{equation} 
\bold{\Lambda} =
 \left(
    \begin{array}{r@{}c|c@{}l}
  &    \begin{smallmatrix}
       1 & & 0 \\
          &\ddots&\\
        0 & & 1 \rule[-1ex]{0pt}{2ex}
      \end{smallmatrix} 
  & \mbox{\huge0} 
  & \rlap{\kern5mm$M$}\\
      \hline
  &    \mbox{\huge0} &  
       \begin{smallmatrix}\rule{0pt}{2ex}
        \lambda & & 0 \\
          &\ddots&\\
        0 & & \lambda
      \end{smallmatrix}    &  \rlap{\kern5mm $N$}
    \end{array} 
\right), \nonumber
\end{equation} 
and $\lambda=\gamma_1/\gamma_0$.
Method of multipliers for the scaled-form Lagrangian reads as
\begin{subequations}
\label{main2}
\begin{eqnarray}
(\bold{m}^{k+1},\bold{u}^{k+1}) &=& \underset{\bold{m},\bold{u}}{\arg\min} ~~ 
 \frac{1}{2}  \|\bold{\Lambda}[\bold{F(m)u}-\bold{s}-\bold{s}^k]\|_2^2, \label{main2a} \\
\bold{s}^{k+1} &=& \bold{s}^{k} +   \bold{s}-\bold{F}(\bold{m}^{k+1})\bold{u}^{k+1}, \label{main2b}
\end{eqnarray} 
\end{subequations}
with  $\bold{s}^0=0$. \\

The scaled-form Lagrangian, equation~\ref{ALS}, is similar to the penalty function of \citet{vanLeeuwen_2016_PMP}, equation~\ref{leeuwen_FWI}, except that the original right-hand side $\bold{s}$ (which gathers $\bold{d}$ and $\bold{b}$ according to our compact notation) is updated at each iteration with the scaled dual variable $\bold{s}^k$. This scaled dual variable gathers the running sum of the constraint violations of previous iterations, namely, the data and source residuals. 
This right-hand side (or residual) updating can be viewed as the counterpart of the residual updating performed at each iteration of the reduced FWI, this residual updating being limited to the data misfit term in this later case (since the wave equation error is 0 in reduced FWI).  \\
\subsubsection*{Alternating-direction method of multipliers}
The joint update of the primal variables ($\bold{m}$,$\bold{u}$) in equation \ref{main2a} is computationally too intensive. 
A splitting method can be useful here to decouple the optimization tasks with respect to $\bold{m}$ and $\bold{u}$ \citep[the readers can refer to ][ for an overview of splitting methods]{Glowinski_2017_SMC}.
Instead of joint updating, the alternating-direction method of multipliers (ADMM) \citep{Boyd_2011_DOS} 
updates the variables $\bold{m}$ and $\bold{u}$ with Gauss-Seidel like iterations,  i.e., fixing one variable and solving for the other.
Using ADMM, we end up with the following iteration:
%
\begin{subequations}
\label{main30}
\begin{eqnarray}
\bold{u}^{k+1} &= & \underset{\bold{u}}{\arg\min} ~~ \frac{1}{2}  \|\bold{\Lambda}[\bold{F}(\bold{m}^k)\bold{u}-\bold{s}-\bold{s}^k]\|_2^2 \label{main3a0} \\
\bold{m}^{k+1} &= & \underset{\bold{m}}{\arg\min} ~~ 
 \frac{1}{2}  \|\bold{\Lambda}[\bold{F(m)u}^{k+1}-\bold{s}-\bold{s}^k]\|_2^2 \label{main3b0} \\
\bold{s}^{k+1} &=& \bold{s}^{k} +  \bold{s} - \bold{F}(\bold{m}^{k+1})\bold{u}^{k+1}.      \label{main3d0}
\end{eqnarray} 
\end{subequations} 
ADMM algorithm has decoupled the full problem into two subproblems associated with primal variables $\bold{u}$ and $\bold{m}$.  Compared to classical ADMM for separable convex problems and linear constraints, the non-separability between wavefield reconstruction and parameter estimation is managed by passing the solution of one primal subproblem as a passive variable for the next subproblem. We will show soon that this alternating-direction approach linearizes the primal subproblem for $\bold{m}$ around the reconstructed wavefield $\bold{u}$. As such, our approach can also be seen as a linearly constrained augmented Lagrangian method, the linearization simply resulting from the splitting procedure \citep[][ pages 522-523]{Nocedal_2006_NOO}. The fact that the primal update of one subproblem is passed to the next subproblem implies obviously that the two subproblems are solved in sequence rather than in parallel as in classical ADMM. \\
For reasons that will become apparent below, we reintroduce explicitly the two components $\bold{d}$ and $\bold{b}$ of the scaled dual vector $\bold{s}$ in the notations
\begin{subequations}
\label{main3}
\begin{eqnarray}
\bold{u}^{k+1} &= &   \underset{\bold{u}}{\arg\min} ~~
 \frac{1}{2}\|\bold{Pu}-\bold{d}-\bold{d}^k\|_2^2+\frac{\lambda}{2}  \|\bold{A}(\bold{m}^{k})\bold{u}-\bold{b}-\bold{b}^k\|_2^2 \label{main3a} \\
\bold{m}^{k+1} &= & \underset{\bold{m}}{\arg\min} ~~ \|\bold{A}(\bold{m})\bold{u}^{k+1}-\bold{b}-\bold{b}^k\|_2^2 \label{main3b} \\
\bold{d}^{k+1} &=& \bold{d}^k + \bold{d}-\bold{Pu}^{k+1} \label{main3d}\\
\bold{b}^{k+1} &=& \bold{b}^k +\bold{b}- \bold{A}(\bold{m}^{k+1})\bold{u}^{k+1}.  \label{main3e}
\end{eqnarray} 
\end{subequations}  
%
A second modification related to classical ADMM resides in the updating of the dual variables.
In classical ADMM, the dual variables are updated only once per iteration after the residuals generated by each primal subproblem have been gathered, equations~\ref{main3d}-\ref{main3e}.    
As reviewed by \citet[][ page 23]{Boyd_2011_DOS} and \citet[][ page 3]{Esser_2009_ALA}, a variant of ADMM, referred to as the Peaceman-Rachford 
splitting method (PRSM) \citep{peaceman_1955_PRA}, consists of updating the Lagrange multipliers several times, once after the update of each primal variable (see \citet[][ equations 1.3 versus 1.4]{He_2014_SPR}.
One issue with PRSM relative to ADMM is to require more restrictive assumptions to ensure its convergence, while it is always faster than ADMM whenever it is convergent \citep{He_2014_SPR}.
This issue prompts \citet{He_2014_SPR} to implement relaxation factors (or, step lengths) $\alpha\in (0,1)$ to guarantee the strict contraction of the PRSM iterative
sequence. 

Accordingly, applying the strictly contractive PRSM algorithm to our minimization problem gives
\begin{subequations}
\label{main4}
\begin{eqnarray}
\bold{u}^{k+1} &= & \underset{\bold{u}}{\arg\min} ~~ 
 \frac{1}{2}\|\bold{Pu}-\bold{d}-\bold{d}^k\|_2^2+\frac{\lambda}{2}  \|\bold{A}(\bold{m}^{k})\bold{u}-\bold{b}-\bold{b}^k\|_2^2 \label{main4a} \\
 \bold{d}^{k+1} &=& \bold{d}^k + \bold{d}-\bold{Pu}^{k+1} \label{main4b}\\
 \bold{b}^{k+\frac{1}{2}} &=& \bold{b}^k + \alpha[\bold{b}-\bold{A}(\bold{m}^{k})\bold{u}^{k+1}]  \label{main4c} \\
\bold{m}^{k+1} &= & \underset{\bold{m}}{\arg\min} ~~
 \|\bold{A}(\bold{m})\bold{u}^{k+1}-\bold{b}-\bold{b}^{k+\frac{1}{2}}\|_2^2 \label{main4d} \\
\bold{b}^{k+1} &=& \bold{b}^{k+\frac{1}{2}} + \alpha[\bold{b}-\bold{A}(\bold{m}^{k+1})\bold{u}^{k+1}]  \label{main4e}
\end{eqnarray} 
\end{subequations}
where $\alpha\in (0,1)$. 

The dual variable $\bold{b}^k$ is updated twice with a step $\alpha$, once after the update of $\bold{u}$ and once after the update of $\bold{m}$, since it depends upon both primal variables, while $\bold{d}^k$ is updated only once after the update of $\bold{u}$. In the following numerical section, we use the PRSM formulation, equation~\ref{main4} at the expense of ADMM, equation~\ref{main3}, as we found that it was leading to a faster convergence provided that $\alpha < 1$ to guarantee stability. 

We now provide the closed-form solution of the two primal subproblems for $\bold{u}$ and $\bold{m}$ and discuss the linearization of the subproblem for $\bold{m}$ resulting from the operator splitting.
\subsubsection{Optimization over $\bold{u}$}
The primal subproblem associated with $\bold{u}$, equations \ref{main3a}-\ref{main4a}, 
is a quadratic optimization problem whose closed-form solution is given by
\begin{equation}  \label{solve_sub_uu} 
\bold{u}^{k+1} ~= ~\Big[\bold{P}^T\bold{P}+\lambda \bold{A}(\bold{m}^k)^T \bold{A}(\bold{m}^k)\Big]^{-1}\Big[ \bold{P}^T [\bold{d}+\bold{d}^k]+\lambda  \bold{A}(\bold{m}^k)^T [\bold{b+b}^k] \Big].
\end{equation}

\subsubsection{Optimization over $\bold{m}$}
The primal subproblem for $\bold{m}$, equations \ref{main3b}-\ref{main4d},  is more challenging due to the nonlinearity of the forward operator $\bold{A}$ in $\bold{m}$.
However, exploiting the special structure of $\bold{A}$, equation~\ref{helmholtz}, it is straightforward to show that the operator splitting has linearized the optimization subproblem for $\bold{m}$ around $\bold{u}$. \\
Indeed, after rewriting $\bold{A(m)u}^{k+1}$ as
\begin{eqnarray}
 \label{nonlinear_op}
\bold{A}(\bold{m})\bold{u}^{k+1} &=& \bold{\Delta}\bold{u}^{k+1} + \omega^2  \bold{B}\bold{C} (\bold{m})\text{diag}(\bold{m}) \bold{u}^{k+1}, \nonumber \\
&\approx & \bold{\Delta} \bold{u}^{k+1} + \bold{L}(\bold{u}^{k+1})\bold{m},
\end{eqnarray}
where
\begin{eqnarray}
\label{eql}
 \bold{L}(\bold{u}^{k+1}) &=& \omega^2  \bold{B}\bold{C} (\bold{m}^k)\text{diag}(\bold{u}^{k+1}),
\end{eqnarray}
we can recast the subproblem \ref{main4d}  as the following linear problem 
\begin{equation} \label{uml}
\bold{m}^{k+1} ~=~ \text{arg} \min_{\bold{m}}  \|\bold{L}(\bold{u}^{k+1}) \bold{m}-\bold{y}^{k}\|_2^2,
\end{equation}
where
\begin{equation}
\bold{y}^{k}=\bold{b}+\bold{b}^{k+\frac12} - \bold{\Delta}\bold{u}^{k+1}.
\end{equation}
Note that we have also linearized the operator $\bold{A}$ with respect to $\bold{m}$ by building the matrix $\bold{C}$ from $\bold{m}^k$ in
equation~\ref{eql} to manage potential nonlinear boundary conditions. However, 
in case of PML absorbing boundary conditions, this extra linearization is not used.
The linear problem, equation~\ref{uml}, has a closed-form solution given by
\begin{equation} \label{mi}
\begin{split}
\bold{m}^{k+1}~= ~\big[\bold{L}(\bold{u}^{k+1})^T \bold{L}(\bold{u}^{k+1})\big]^{-1} \big[\bold{L}(\bold{u}^{k+1})^T \bold{y}^{k}\big].
\end{split}
\end{equation}
Note that the Gauss-Newton Hessian ($\bold{L}^T \bold{L}$) in equation~\ref{mi} is diagonal if the mass term is not spread over the coefficients of the stencil, which might be not necessary to solve the inverse problem, while beneficial for the forward problem \citep[][ their equation 8]{VanLeeuwen_2013_MLM}. \\
Now that we have linearized the minimization problem for $\bold{m}$ around $\bold{u}$ through the splitting procedure, the two subproblems embedded in the PRSM method can be discussed in the framework of linear inverse problem theory. The right-hand side correction terms $\bold{b}^k$ and $\bold{d}^k$ in the objective, equation~\ref{main4}, gather the running sum of the constraint violations (data and source residuals) of previous iterations. For linear problems, these accumulated residuals are linearly related to the running sum of the solution refinements (model perturbations in the classical FWI terminology) performed at previous iterations. Therefore, these right-hand side correction terms gather in the objective the imprint of the solution refinements performed at previous iterations such that the solution refinement performed at the current iteration relies only on the residual errors (see equation~\ref{eqa10} in the appendix). 
This iterative defect correction leads to the error forgetting property discussed by \citet{Yin_2013_EFG} in the frame of Bregman iterations, which means that the error correction performed at the current iteration is made independent to the error corrections performed at previous iterations.
Some applications of iterative error correction in geophysics are for example least-squares migration \citep[e.g.][]{Lambare_1992_IAI,Jin_1992_TDA,Ribodetti_2011_JRB}, which can be viewed as the linear counterpart of the iterative data residual minimization performed by nonlinear FWI. Another application is presented by \citet{Gholami_2017_CNA} who perform  nonlinear amplitude versus offset (AVO) inversion, using a linearized Zoeppritz equation to formulate the iterative linearized inverse problem, while computing the errors with the nonlinear equation. The principles of the iterative solution refinement for solving linear inverse problems are reviewed in the appendix. \\
Based on the above, the proposed PRSM algorithm for FWI is summarized in Algorithm~\ref{Alg2cont}. Note that, in the outer loop (lines 3-12), the algorithm performs one iteration of the wavefield reconstruction (lines 4-6) and one iteration of parameter estimation (lines 9-10), while several inner iterations of each subproblem can be viewed (these inner iterations are omitted in Algorithm~\ref{Alg2cont} for sake of compactness).
The usefulness of these inner iterations is discussed in the next section. \\
In the next section, we will also consider that the results obtained with IR-WRI when right-hand side updating is not activated are representative of those that would be obtai, is identically zeroned with WRI as shown by {equations~\ref{ALS}-\ref{main2}},

\SetAlgorithmName{Algorithm 1}{}{}
\begin{algorithm}[H]
\label{Alg2cont}
\small
\setstretch{1.35}
Initialize: set $k=0$, $\bold{d}^{k} = \bold{0}$, $\bold{b}^{k} = \bold{0}$.\\
Input: $\bold{m}^0$\\
\While{convergence criteria  not satisfied}{
$\bold{u}^{k+1} \leftarrow  \text{update according to equation}~~ \ref{solve_sub_uu}$\\
$\bold{d}^{k+1} \leftarrow  \bold{d}^k + \bold{d} - \bold{Pu}^{k+1} $\\
$\bold{b}^{k+\frac{1}{2}} \leftarrow  \bold{b}^k + \alpha\big[\bold{b}-\bold{A}(\bold{m}^{k})\bold{u}^{k+1}\big]$\\
$\bold{y}^{k} \leftarrow \bold{b} + \bold{b}^{k+\frac{1}{2}} - \bold{\Delta}\bold{u}^{k+1}$\\
$\bold{L}(\bold{u}^{k+1}) \leftarrow  \omega^2\bold{B}\bold{C}(\bold{m}^{k})\text{diag}(\bold{u}^{k+1})$\\
$\bold{m}^{k+1} \leftarrow  \text{update according to equation}~~ \ref{mi}$\\
$\bold{b}^{k+1} \leftarrow  \bold{b}^{k+\frac{1}{2}} + \alpha\big[\bold{b}-\bold{A}(\bold{m}^{k+1})\bold{u}^{k+1}\big]$\\
$k \leftarrow  k+1$
 }
 \caption{
FWI algorithm based on the Peaceman-Rachford splitting algorithm.}
\end{algorithm}

\section{Numerical examples}
We compare the performance of WRI and IR-WRI against three 2D mono-parameter synthetic examples: a toy example built with a simple inclusion model 
to compare the convergence history of the two methods and two complex synthetic experiments where the synthetics computed in the starting model are strongly cycle skipped relative to the observables. 
%
For each example, no regularization is used. We only apply bound constraints using the true values of the minimum and maximum velocities as bounds for the velocity-model update.  Because the interfacing of regularization and additional constraints in our method is not the scope of this study, we have not reviewed the bound-constraint implementation in the method section for sake of clarity. We implement bound constraints through an auxiliary variable and add an extra constraint to the problem following the slit Bregman method \citep{Goldstein_2009_SBM}.
The interested reader is also referred to \citep{Maharramov_2015_TVM}, where bound constraints and total-variation regularization are implemented in FWI with the split Bregman method \citep{Goldstein_2009_SBM}.\\  
For all numerical examples the forward modeling is performed with a 9-point stencil finite-difference method implemented with anti-lumped mass and PML absorbing boundary conditions to solve the Helmholtz equation, where the stencil coefficients are optimized to the frequency \citep{Chen_2013_OFD}. In this setting the diagonal matrix $\bold{C}$ contains the damping PML coefficients and does not depend on $\bold{m}$.  With this setting, the update of $\bold{m}$ can be recast as a linear inverse problem, without additional linearization to that induced by the operator splitting (see equation~\ref{mi}). \\
%
%
%
%
%
\subsection{Inclusion model}
\begin{figure}
\begin{center}
\includegraphics[scale=0.6]{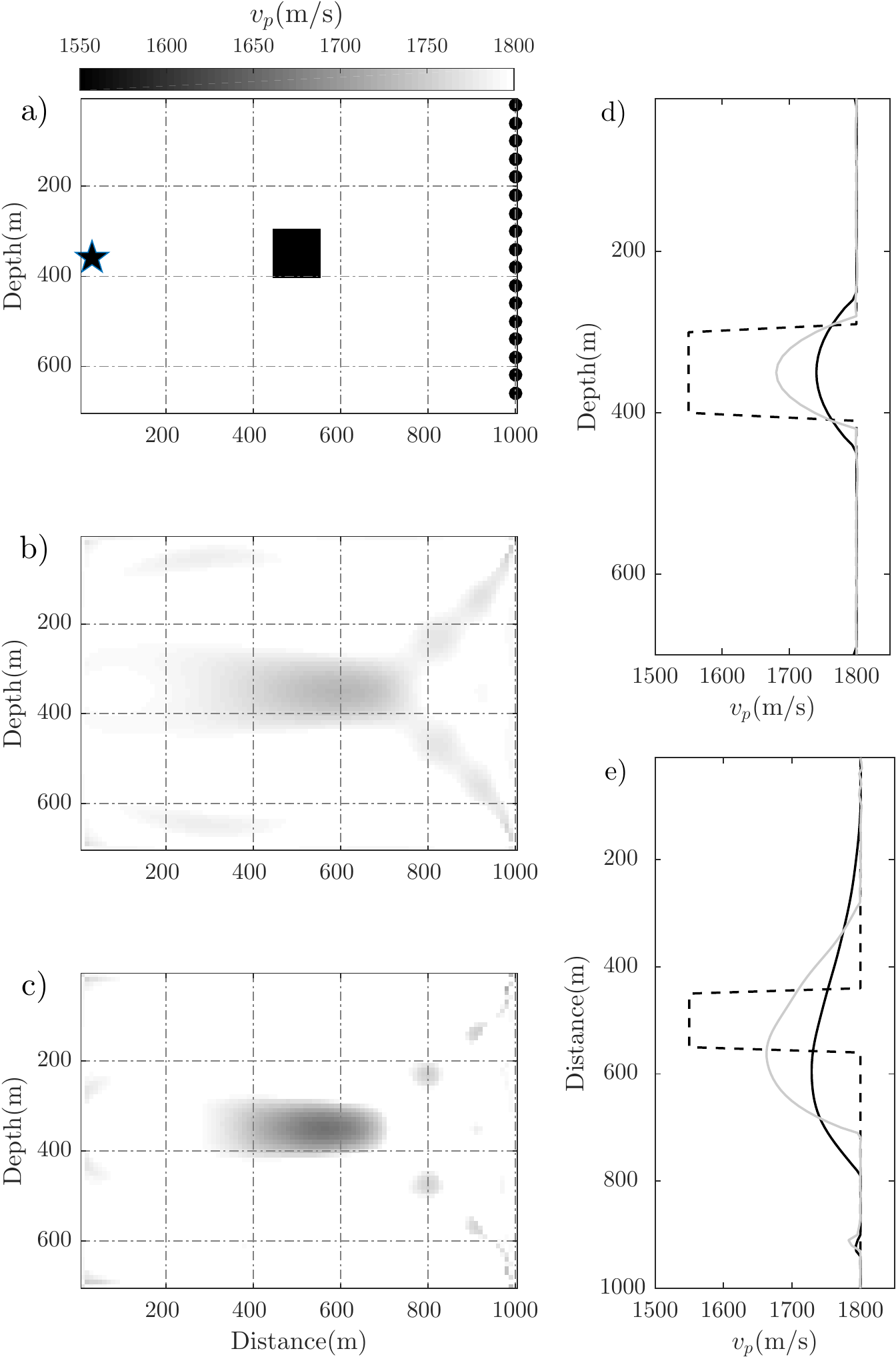}
\caption{Box-shape anomaly example. (a) True model. The star points the source position, while circles denote receivers. (b-c) WRI (b) and IR-WRI (c) models after 100 iterations. (d-e) Direct comparison between velocity models shown in (a-c) along the vertical (d) and horizontal (e) profiles cross-cutting the center of the anomaly. The dashed, black and gray lines correspond to the models shown in (a), (b) and (c), respectively. }
\label{fig:fig_box_new}
\end{center}
\end{figure}

The subsurface model contains a box-shape anomaly of side 100~m embedded in a homogeneous background (Figure~\ref{fig:fig_box_new}a). The dimension of the model are 1000~m in distance and 700~m in depth and the grid spacing is 10~m. The cross-hole acquisition consists of one source at (x,z)=(0~m,350~m) and 18 receivers deployed vertically on the opposite side of the grid (Figure~\ref{fig:fig_box_new}a). The source signature is a Ricker wavelet with a 5~Hz dominant frequency. We start the inversion from the homogeneous background model ($v_p=1800$~m/s) and invert simultaneously three frequency components (2.5, 5 and 7 Hz) with noiseless data. \\
To gain first hints on the relative convergence speed of WRI and IR-WRI, we set the maximum number of iteration to 100 as stopping criterion of iteration (Figure~\ref{fig:fig_box_new}b-c). The value of the penalty parameter $\lambda$, equation \ref{solve_sub_uu}, has been chosen to be a small fraction (1e-4) of the largest eigenvalue of $\bold{A}^{-T}\bold{P}^{T}\bold{P}\bold{A}^{-1}$ according to the scaling proposed by \citep{vanLeeuwen_2016_PMP}. In the following, we denote this maximum eigenvalue by $\mu_1$ according to the notation of \citep{vanLeeuwen_2016_PMP}. The WRI model is clearly smoother than the one obtained with IR-WRI suggesting that iterative solution refinement performed by IR-WRI contributes to decrease the objective faster.
\subsubsection{Sensitivity of WRI and IR-WRI to $\lambda$ and convergence analysis}
\begin{figure}
\begin{center}
\includegraphics[scale=0.65]{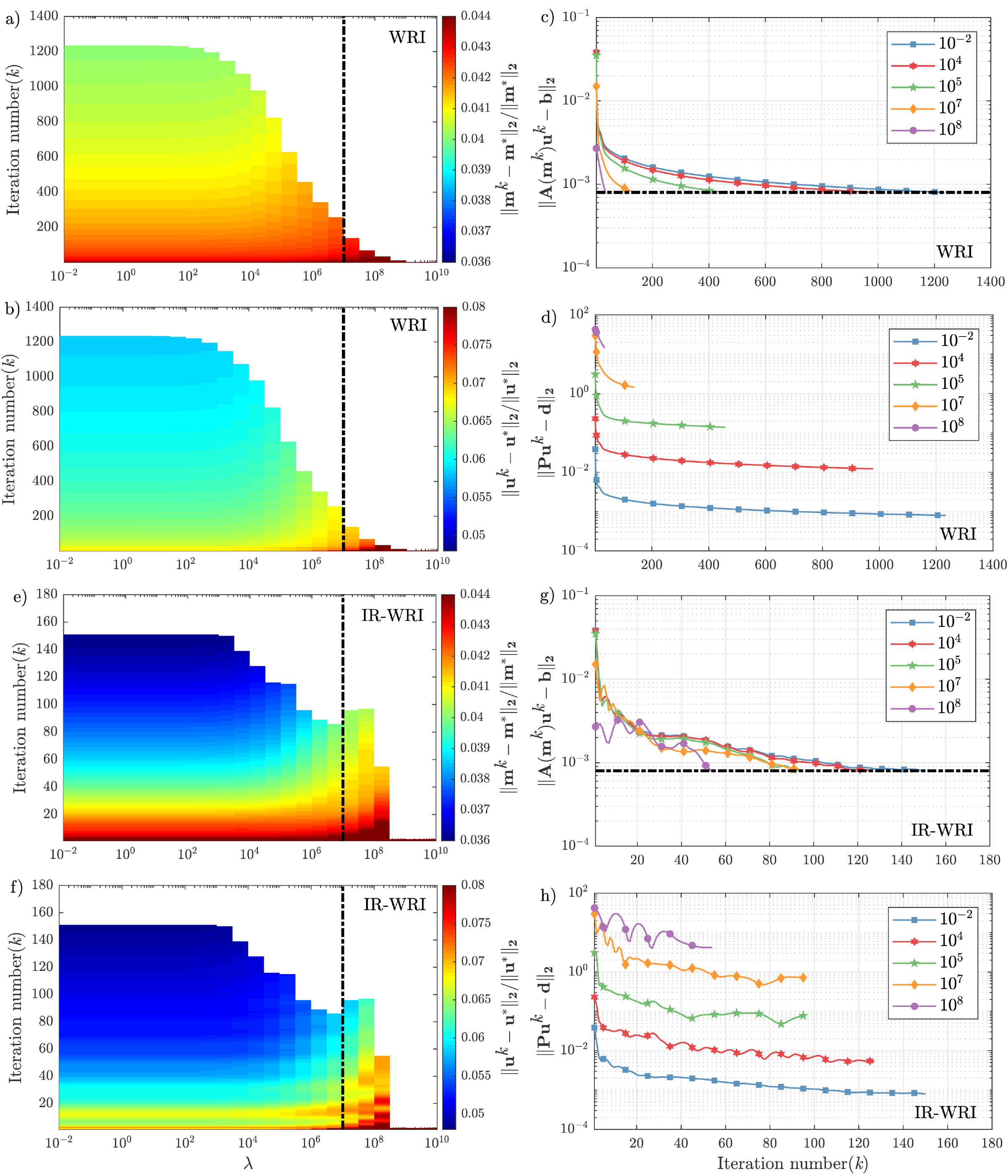}
\caption{Box-shape anomaly example. Sensitivity of WRI and IR-WRI to $\lambda$. (a-d) WRI results. (a-b) Convergence history of (a) model error function $||\bold{m}^k-\bold{m}^*||_2/||\bold{m}^*||_2$ and  (b) wavefield error function $||\bold{u}^k-\bold{u}^*||_2/||\bold{u}^*||_2$ in the ($\lambda$-$k$) plane. (c-d) Convergence history of (c) wave-equation misfit function $||\bold{A(m}^k)\bold{u}^k-\bold{b}||_2$ and (d) data misfit function $||\bold{Pu}^k-\bold{d}||_2$ as function of $k$ for several $\lambda$. (e-h) Same as (a-d) for IR-WRI. The vertical dashed line in (a), (b), (e) and (f) points the
value of the highest eigenvalue of matrix $\bold{A}^{-T}\bold{P}^{T}\bold{P}\bold{A}^{-1}$ (namely, 1e7) and the horizontal dashed line in (c) and (g) shows the preset stopping criteria based upon the wave-equation misfit function $||\bold{A(m}^k)\bold{u}^k-\bold{b}||_2$.}
\label{fig:fig_misfit_box_31}
\end{center}
\end{figure}

\begin{figure}
\begin{center}
\includegraphics[scale=0.7]{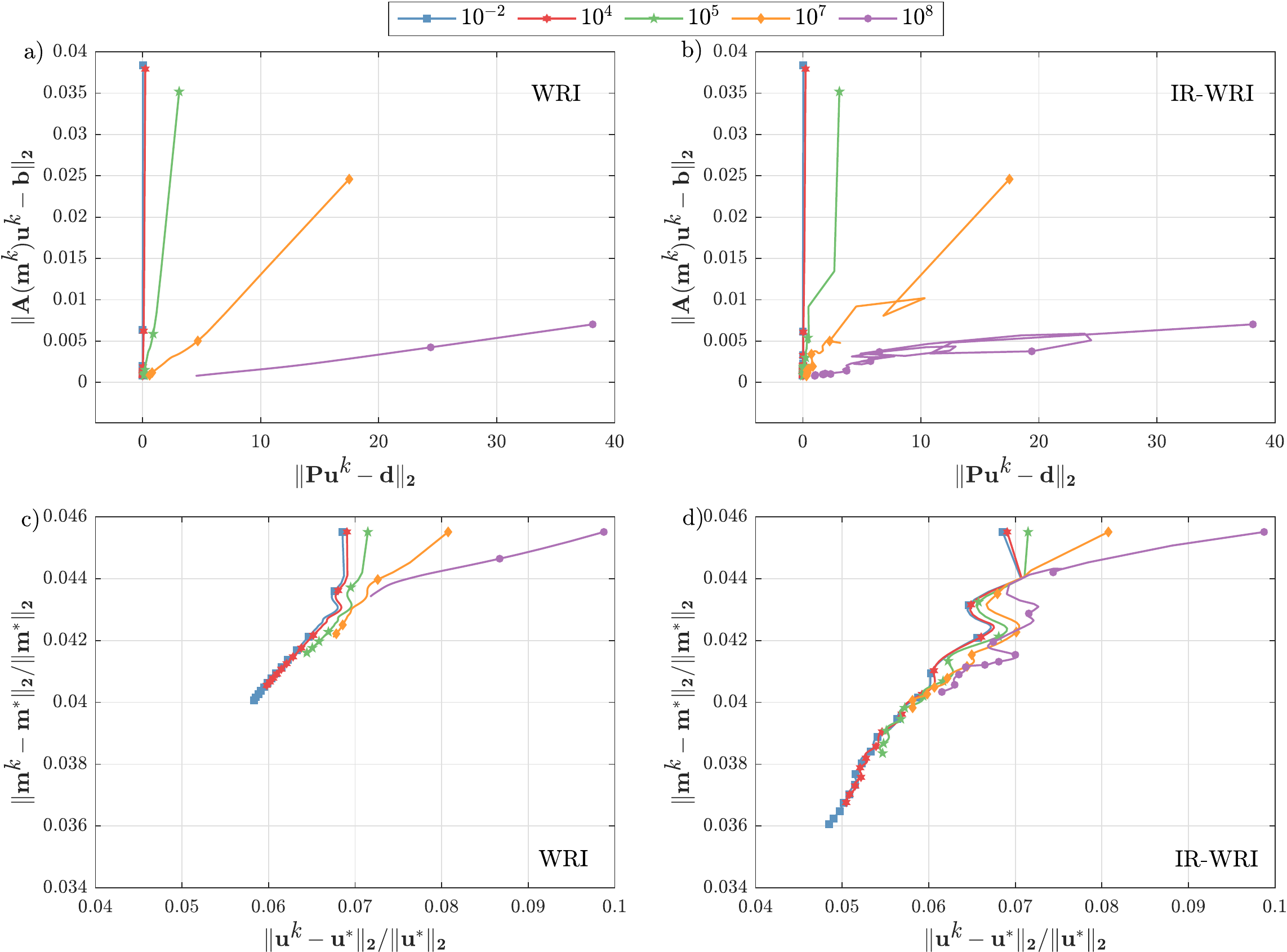}
\caption{Box-shape anomaly example. Convergence history of WRI (a,c) and IR-WRI (b,d) (see also Figure \ref{fig:fig_misfit_box_31}). (a-b) Convergence path in the ($||\bold{Pu}^k-\bold{d}||_2 - ||\bold{A(m}^k)\bold{u}^k-\bold{b}||_2$) plane for several values of $\lambda$. (c-d) Convergence path in the ($||\bold{u}^k-\bold{u}^*||_2/||\bold{u}^*||_2-||\bold{m}^k-\bold{m}^*||_2/||\bold{m}^*||_2$) plane for several values of $\lambda$. The circles denote the iteration count in (a-d).}
\label{fig:fig_Box_res1_res2_plot1}
\end{center}
\end{figure}

\begin{figure}
\begin{center}
\includegraphics[scale=0.6]{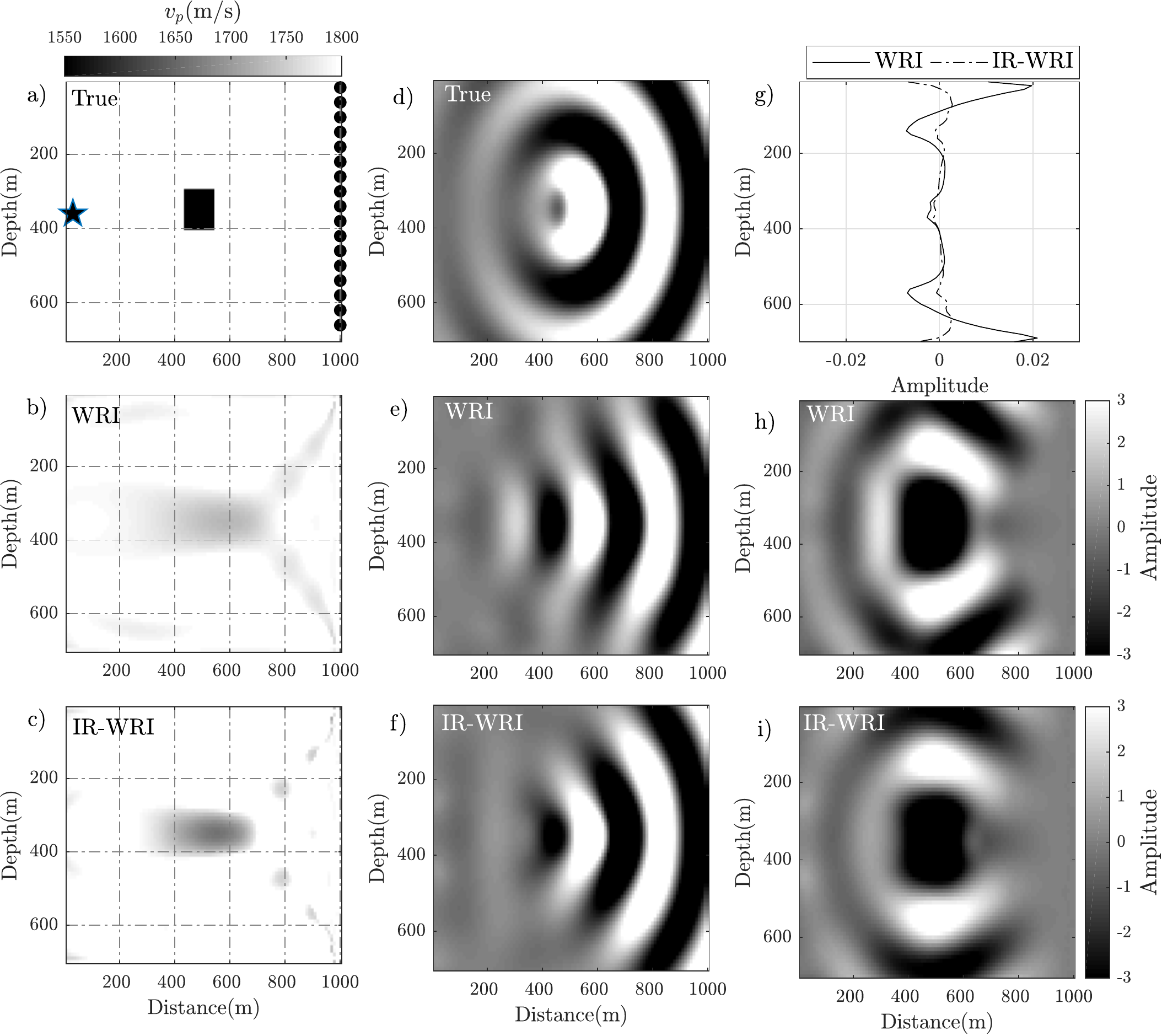}
\caption{Box-shape anomaly example. WRI and IR-WRI final results for $\lambda$=1000 (see Figure \ref{fig:fig_misfit_box_31}a,e). (a) True model. (b-c) Final WRI (b) and IR-WRI (c) models. (d) Wavefield (real part) scattered by the true inclusion for the 7~Hz frequency. (e-f) Wavefield (real part) scattered by the inclusion reconstructed by (e) WRI and (f) IR-WRI for the 7~Hz frequency. (g) Final 7~Hz data residuals $\bold{Pu}^k-\bold{d}$ (real part) at the receiver positions for WRI (solid line) and IR-WRI (dashed line). (h) Difference between scattered wavefields shown in (d) and (e). (i) Same as (h) for scattered wavefields shown in (d) and (f).}
\label{fig:fig_box_with_wavefield2}
\end{center}
\end{figure}

A more comprehensive convergence analysis can be performed by applying WRI and IR-WRI with several values of $\lambda$ ranging between 1e-9 $\mu_1$ and 1e3 $\mu_1$ (Figure~\ref{fig:fig_misfit_box_31}). The value of $\mu_1$ (=1e7) is indicated by the vertical dash line in Figure~\ref{fig:fig_misfit_box_31}a,b,e,f. We stop iteration when a preset value of the wave-equation misfit $\|\bold{A}(\bold{m}^k)\bold{u}^k-\bold{b}\|_2$ is reached  (Figure~\ref{fig:fig_misfit_box_31}c,g, horizontal dash line). 
We compare the relative model error $\|\bold{m}^k - \bold{m}^*\|_2 / \|\bold{m}^*\|_2$ and the relative wavefield error $\|\bold{u}^k - \bold{u}^*\|_2 / \|\bold{u}^*\|_2$ achieved by WRI and IR-WRI as functions of the iteration count $k$ and  $\lambda$ in Figure~\ref{fig:fig_misfit_box_31}. In parallel with this, we compare the convergence history of the wave-equation objective $\|\bold{A}(\bold{m}^k)\bold{u}^k-\bold{b}\|_2$ and the data misfit objective $\|\bold{P}\bold{u}^k-\bold{d}\|_2$ achieved by WRI and IR-WRI as functions of $k$ for several $\lambda$ in Figure~\ref{fig:fig_misfit_box_31}. The true model and the wavefield computed in the true model (namely, the global minimizers) are denoted by $\bold{m}^*$ and $\bold{u}^*$, respectively. \\
%
For small values of $\lambda$ ranging between 1e-9 $\mu_1$ and 1e-4 $\mu_1$, IR-WRI reaches the stopping criterion of iteration after a much smaller number of iterations (typically, one order of magnitude smaller) than WRI (Figure~\ref{fig:fig_misfit_box_31}c,g). The most-accurate minimizers (wavefield and subsurface model) are obtained for these small values of $\lambda$, whatever WRI or IR-WRI is used, because the inversion has to decrease the data residuals significantly before being able to minimize the source residuals (due to the small  weight assigned to the wave equation objective) and fulfill the stopping criterion of iteration accordingly (Figure~\ref{fig:fig_misfit_box_31}a,b,e,f). However, the minimizers estimated by IR-WRI are significantly more accurate than those of WRI for a fixed value of $\lambda$ (compare Figure~\ref{fig:fig_misfit_box_31}a and ~\ref{fig:fig_misfit_box_31}e, Figure~\ref{fig:fig_misfit_box_31}b and \ref{fig:fig_misfit_box_31}f).
This highlights that IR-WRI has decreased data residuals not only faster but also more significantly than WRI before reaching the stopping criterion of iteration.  \\
This more efficient data residual minimization performed by IR-WRI likely results from the joint updating of the data and source residuals in iterations. However, one may wonder whether this joint residual updating performed by IR-WRI presents the risk that the wave-equation based stopping criterion of iteration is prematurely satisfied before achieving a sufficient data fit. The numerical results suggest that this does not to occur because each time the wave equation error is decreased, this pushes back the optimization against the data fitting constraint, hence re-balancing the relative weight of the two objectives. This interpretation highlights the self-adaptivity of the augmented Lagrangian method to process competing constraints by iteratively updating the Lagrange multipliers (see \citep[][ Section 12.8]{Nocedal_2006_NOO} for the intuitive significance of Lagrange multipliers). As such, the augmented Lagrangian method written in a scaled form, equation~\ref{ALS}, can be viewed as a self-adaptive penalty method, where the iterative updating of the Lagrange multipliers provides a systematic and easy-to-implement alternative to an ad hoc penalty parameter continuation strategy. \\
When $\lambda$ is decreased, the data fit achieved by WRI and IR-WRI becomes closer and closer, since the small weight assigned to the wave-equation objective leaves some time for WRI and IR-WRI to get close to a local minimum where the objective is nearly flat before satisfying the stopping criterion of iteration. At the same time, the number of iterations performed by the two approaches becomes increasingly different because the poorly-conditioned WRI needs a prohibitively-high number of iterations before the stopping criterion of iteration is fulfilled. \\
%
%
When higher values of $\lambda$ are used up to a value of 1e3 $\mu_1$, WRI and IR-WRI satisfy the stopping criterion after a smaller number of iterations (down to around 60 iterations for $\lambda$ = 10 $\mu_1$) (Figure~\ref{fig:fig_misfit_box_31}c,g) before an acceptable data fit has been achieved (Figure~\ref{fig:fig_misfit_box_31}d,h). This degraded data fit translates into minimizers of more limited accuracy (Figure~\ref{fig:fig_misfit_box_31}a,b,e,f). However, this degradation of the data fit and minimizer quality is less significant in IR-WRI than in WRI due to the above-mentioned self-adaptive management of the two objectives performed by IR-WRI.   \\
It is also instructive to look at the joint evolution in iterations of the data misfit and the wave-equation error on the one hand (Figure~\ref{fig:fig_Box_res1_res2_plot1}(a-b), and the wavefield and subsurface model errors on the other hand (Figure~\ref{fig:fig_Box_res1_res2_plot1}(c-d). Figure~\ref{fig:fig_Box_res1_res2_plot1}(a-b) shows that small values of $\lambda$ allow for significant wave-equation error during early iterations, which is consistent with the governing idea of expanding the search space.
%
The more complex zigzag path followed by IR-WRI relative to WRI in the ($||\bold{Pu}^k-\bold{d}||_2 - ||\bold{A(m}^k)\bold{u}^k-\bold{b}||_2$) plane for intermediate values of $\lambda$ (this trend is also illustrated in Figure~\ref{fig:fig_misfit_box_31}c,d,g,h by the non monotonicity of the IR-WRI convergence curves) highlights how the joint updating of the data misfit and wave-equation error dynamically balances the weight of the two objectives in iterations. This zigzag convergence trend translates also into more complex path in the $(||\bold{u}^k-\bold{u}^*||_2/||\bold{u}^*||_2-||\bold{m}^k-\bold{m}^*||_2/||\bold{m}^*||_2)$ plane, which highlights how the solution refinement is pushed toward the wavefield reconstruction or the velocity model estimation according to the self-adapting weighting of the two objectives (Figure~\ref{fig:fig_Box_res1_res2_plot1}(c-d)). Figure~\ref{fig:fig_Box_res1_res2_plot1}(c-d) also shows that the accuracy gap between the WRI and IR-WRI minimizers becomes increasingly significant as $\lambda$ is decreased, while the data fit achieved by the two methods becomes closer (Figure~\ref{fig:fig_Box_res1_res2_plot1}(a-b)). This highlights the ability of IR-WRI to keep on refining significantly the wavefield and velocity model when the optimization becomes close to a local minimum where the objective function is increasingly flat. \\
Finally, we show the final wavefields and velocity models as well as the final data misfit obtained with WRI and IR-WRI for $\lambda$=1e-4 $\mu_1$ (Figure~\ref{fig:fig_box_with_wavefield2}).
It is reminded that the final wave-equation error is the same for WRI and IR-WRI as it has been used as stopping criterion of iteration, and hence is not shown here. Moreover, $\lambda$=1e-4 $\mu_1$ provides the best-trade off between the number of iterations and the solution accuracy (Figure~\ref{fig:fig_misfit_box_31}a,e).
WRI fails to  build a velocity model as well resolved as the one found by IR-WRI (Figure~\ref{fig:fig_box_with_wavefield2}(b-c)), as when the maximum number of iteration was limited to 100 (Figure~\ref{fig:fig_box_new}).  Increasing $\lambda$ would not solve the problem as it would lead to a poor data fit. Here, the WRI failure does not obviously result from cycle skipping as this case study lies in the linear regime of classical FWI but instead highlights the inability of the penalty method to retrieve accurate minimizer with a fixed $\lambda$. Instead, the iterative updating of the constraint violations achieved by the augmented Lagrangian method provides an efficient and automatic substitute to the adaptive tuning of the penalty parameter to jointly satisfy the two competing data-fitting and wave-equation objectives at the final iteration.
The final WRI and IR-WRI  (scattered) wavefield are shown in Figure~\ref{fig:fig_box_with_wavefield2}(e-f). Due to the transmission acquisition geometry, only the forward-scattered wavefield is reconstructed successfully. Figure~\ref{fig:fig_box_with_wavefield2}(g-i) confirms that IR-WRI outperforms WRI to fit both wavefield and data.
%
%
\subsubsection{ADMM versus Peaceman-Rashford dual updating}
\begin{figure}
\begin{center}
\includegraphics[scale=0.7]{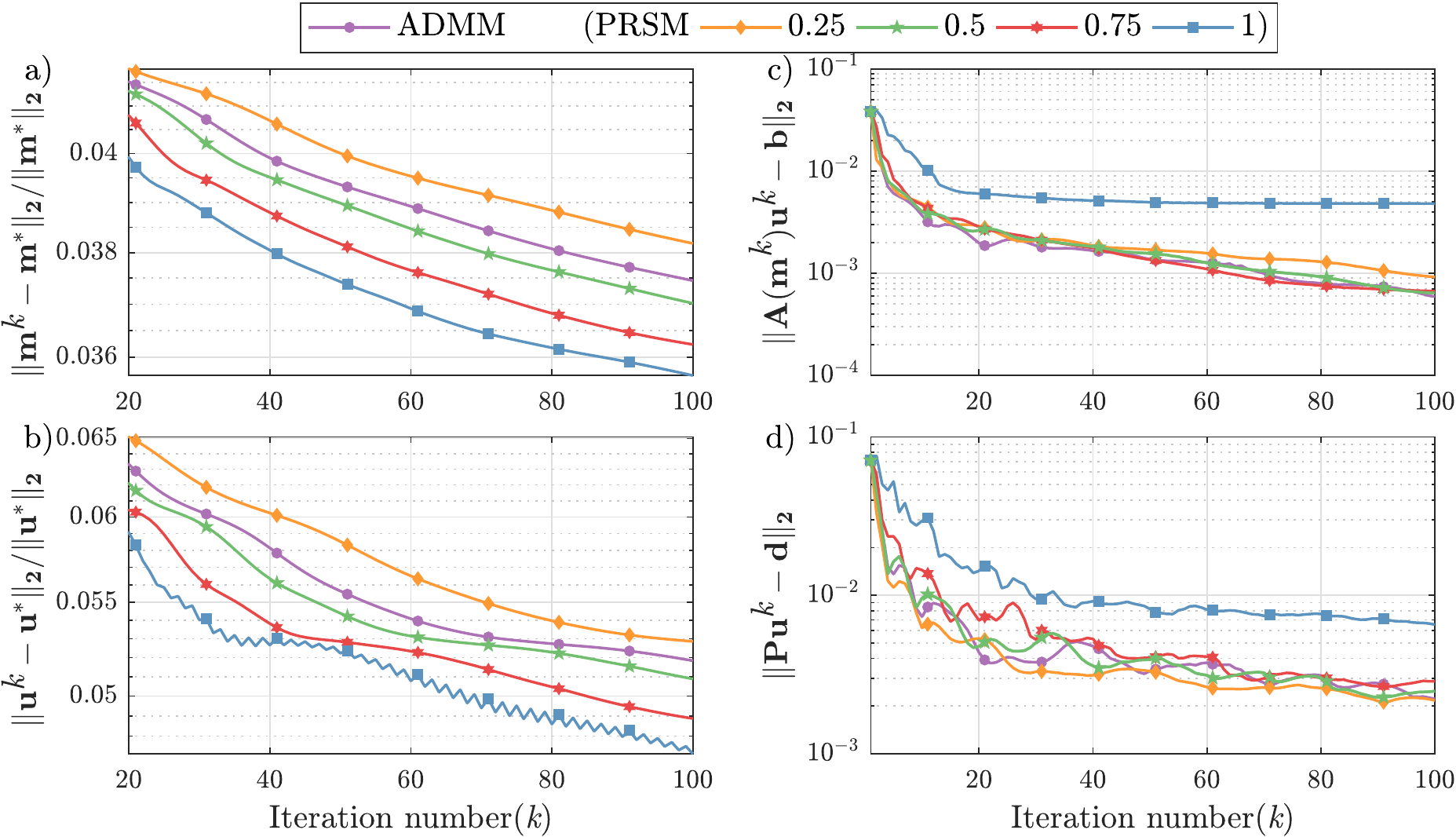}
\caption{Box-shape anomaly example. Comparison between ADMM and PRSM and sensitivity of PRSM to the step length $\alpha$ applied to the $\bold{b}$ updates. (a-d) Convergence history of (a) the model error function $||\bold{m}^k-\bold{m}^*||_2/||\bold{m}^*||_2$, (b) the wavefield error function $||\bold{u}^k-\bold{u}^*||_2/||\bold{u}^*||_2$, (c) the wave-equation  misfit function $||\bold{A(m}^k)\bold{u}^k-\bold{b}||_2$, and (d) the data misfit function $||\bold{Pu}^k-\bold{d}||_2$.}
\label{fig:fig_Box_weight1}
\end{center}
\end{figure}

In IR-WRI algorithm based on the Peaceman-Rachford splitting (see section {\it{Alternating-direction method of Lagrange multipliers}}), the dual variable $\bold{b}^k$ is updated twice with a step length equal to $\alpha$, equation~\ref{main4}. To show the relevance of this double dual updating, we perform IR-WRI for several weights $\alpha$ and plot the model and wavefield errors as well as the wave-equation and data misfit as functions of $k$ in Figure~\ref{fig:fig_Box_weight1}. We also compare the results of PRSM with ADMM, this later consisting of updating $\bold{b}^k$ only once with a weight equal to 1, equation~\ref{main3}. 
The PRSM results obtained with $\alpha$=1 can lead to an improved minimizers (Figure~\ref{fig:fig_Box_weight1}(a-b)), although this setting leads potentially to instabilities. These instabilities are illustrated by a poorer data fit and higher wave equation error when $\alpha$=1 in Figure~\ref{fig:fig_Box_weight1}(c-d)). Overall, PRSM builds improved minimizers when $\alpha$ becomes closer to 1. Also, the PRSM minimizers obtained with $\alpha \geq 0.5$ are more accurate than those obtained with ADMM. Accordingly, we perform the following Marmousi II and BP salt experiments with PRSM using $\alpha$=0.5.

%
\subsubsection{Waveform inversion linearization and iteration strategy}
\begin{figure}
\begin{center}
\includegraphics[scale=0.7]{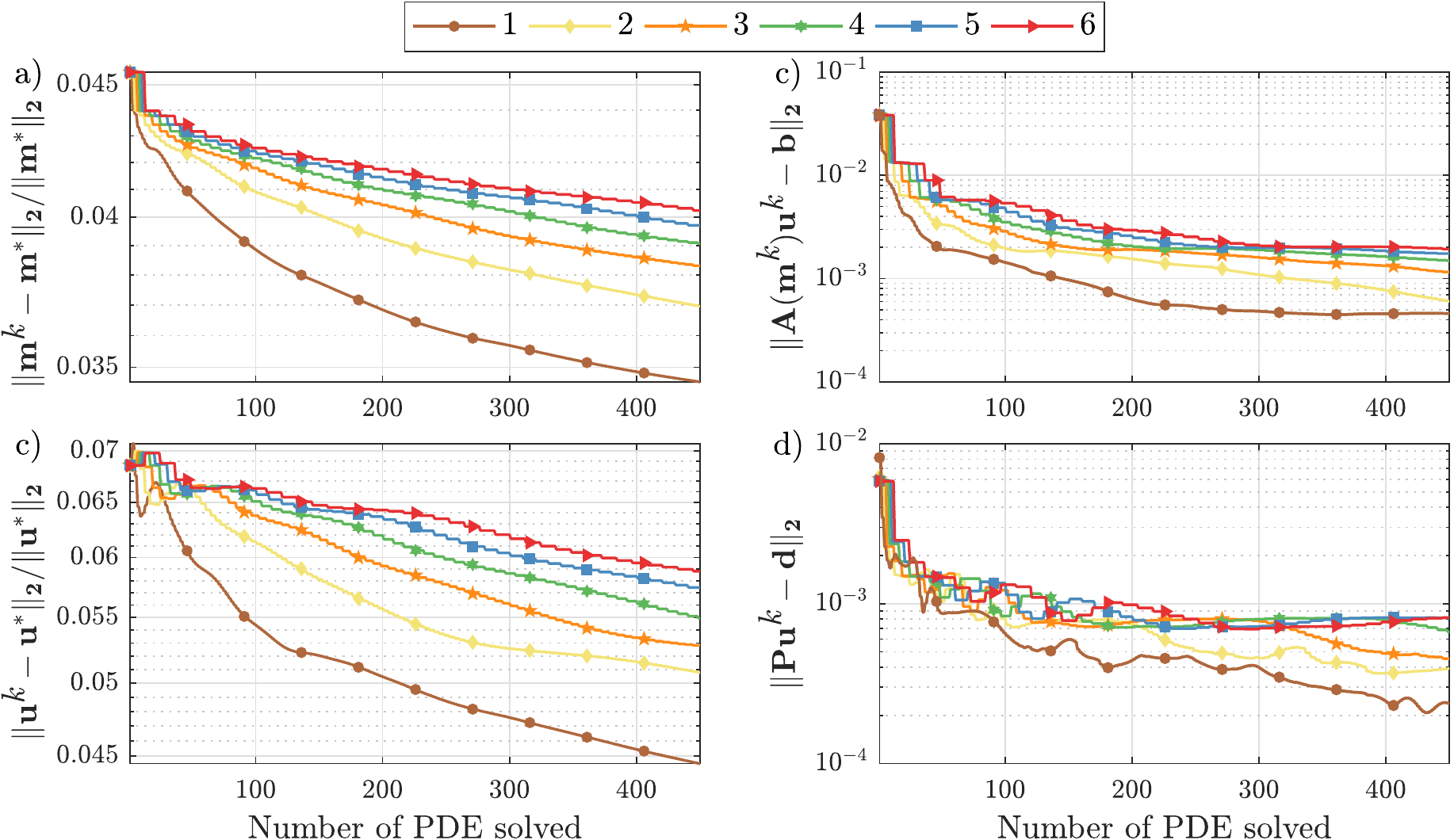}
\caption{Box-shape anomaly example. Sensitivity of the WRI convergence to the number $n$ of inner iterations (see text for details). (a-b) Convergence history of (a) the model error function $||\bold{m}^k-\bold{m}^*||_2/||\bold{m}^*||_2$, (b) the wavefield error function $||\bold{u}^k-\bold{u}^*||_2/||\bold{u}^*||_2$, (c) the wave-equation  misfit function $||\bold{A(m}^k)\bold{u}^k-\bold{b}||_2$, and (d) the data misfit function $||\bold{Pu}^k-\bold{d}||_2$.  The horizontal axis is labeled with the number of PDE solve. In the figure, $n$ ranges from 1 to 6 and the same $n$ is used during the update of $\bold{u}$ and $\bold{m}$.}
\label{fig:fig_innerloop2_new1}
\end{center}
\end{figure}

%
The IR-WRI algorithm potentially embeds two levels of nested iterations, Algorithm~\ref{Alg2cont}. The outer loop over $k$ manages the iterations of the nonlinear multi-variate optimization. Inside one cycle, we can refine the wavefield several times through the iterative updating of $\bold{d}$ and $\bold{b}$ for the current $\bold{m}^k$, lines 4-6 in Algorithm~\ref{Alg2cont}, before refining the parameters several times through the iterative updating of $\bold{b}$ for the current $\bold{u}^{k+1}$, lines 9-10 in Algorithm~\ref{Alg2cont}. For sake of compactness, these two successive inner loops are not explicitly written in Algorithm~ \ref{Alg2cont}, where one iteration of the wavefield reconstruction and parameter estimation are performed per cycle. Let's denote by $n$ the number of inner iterations for the wavefield reconstruction and the parameter update. We perform IR-WRI for several values of $n$ and show the wavefield and model errors as well as the wave-equation misfit and data misfit as functions of the number of PDE resolution (each inner iteration pair requires one additional PDE resolution) in Figure~\ref{fig:fig_innerloop2_new1}. The results clearly show that the best choice is $n$=1. This can be intuitively understood by reminding that the original nonlinear multivariate optimization for $\bold{u}$ and $\bold{m}$ has been recast as two subproblems for $\bold{u}$ and $\bold{m}$ that are solved in sequence. Moreover, this operator splitting allows for the linearization of the model update around the reconstructed wavefield, equation~\ref{mi}. This operator splitting strategy implies that, during the resolution of each subproblem, the optimization variable is updated from an inaccurate passive variable. Therefore, the resulting residuals cannot be decreased efficiently by inner iterations due to the inaccuracy of the passive variable which is kept fixed. A more efficient procedure is therefore to pass the residuals and the optimization variable of one subproblem to the next one as soon as they have been updated by one inner iteration. This conclusion is similar to the one reached by \citep[][section 3.2]{Goldstein_2009_SBM} in the different framework of $l1$-regularized constrained optimization problem solved with the split Bregman method, namely an operator splitting method similar to ADMM.
%
%
%
%
%
%
%
%
%
%
%
%
%
%
%
%
%
%
%
%
%
%
%

\subsection{Marmousi model}
\begin{figure}
\begin{center}
\includegraphics[scale=0.75,trim=0cm 5cm 0cm 0cm]{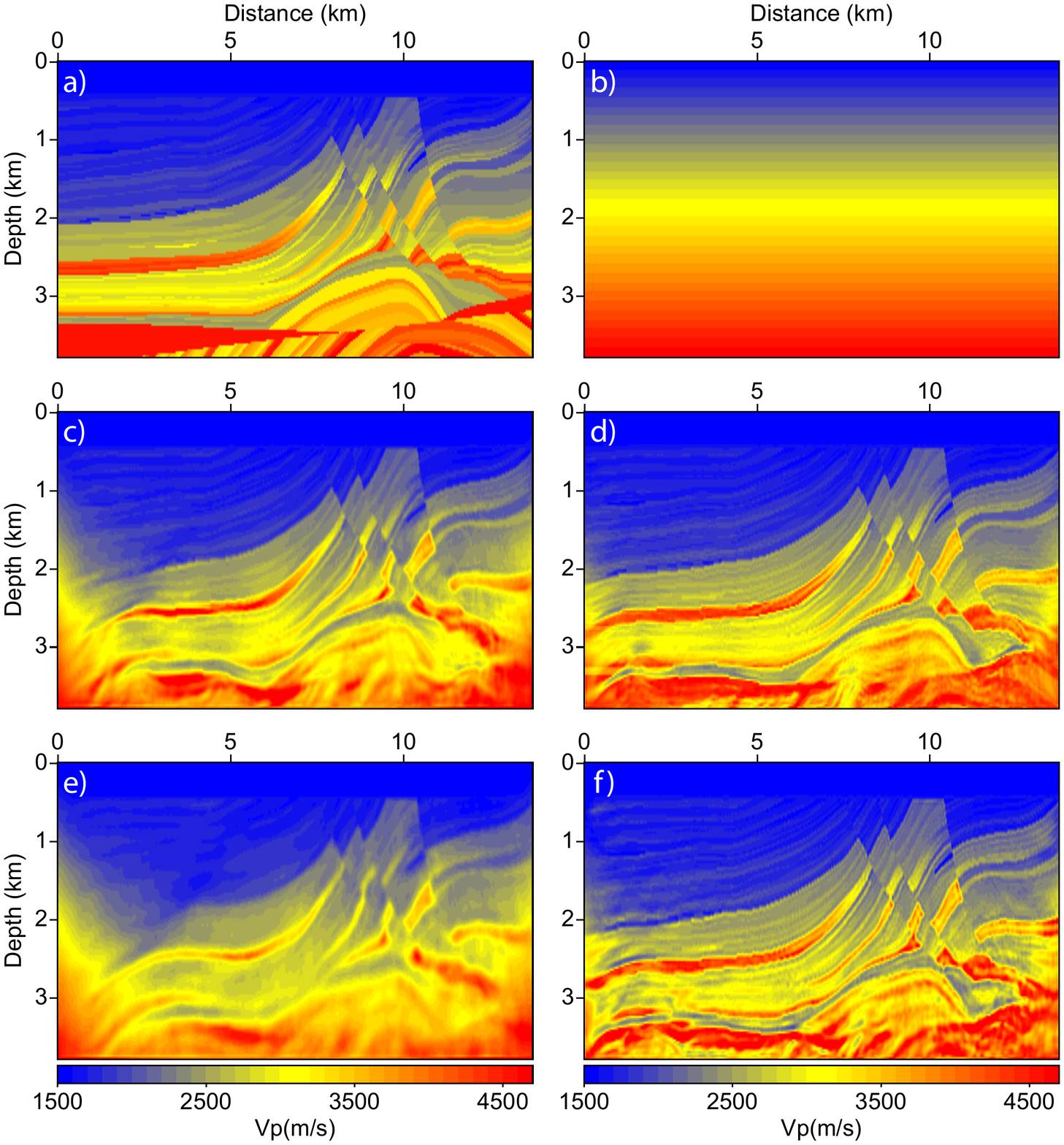}
\caption{Marmousi II example. (a) True Marmousi II model. (b) Initial velocity model. (c-d) For noiseless data, velocity models built by (c) WRI and (d) IR-WRI. (e-f) Same as (c-d) for noisy data with a SNR of 10~db.}
\label{fig:fig_marmousi_full_label}
\end{center}
\end{figure}

We consider now the more complex Marmousi II model, which covers a $13750~m \times 3800~m$ spatial domain (Figure \ref{fig:fig_marmousi_full_label}a). We perform WRI and IR-WRI in the 3~Hz - 15~Hz frequency band and we re-sample the original model on a 25~m grid accordingly. The fixed-spread surface acquisition consists of 137 sources spaced 100~m apart and 548 receivers spaced 25~m apart at the surface.
%
%
%
%
\begin{figure}
\begin{center}
\includegraphics[scale=0.7,trim=0cm 17cm 0cm 0cm]{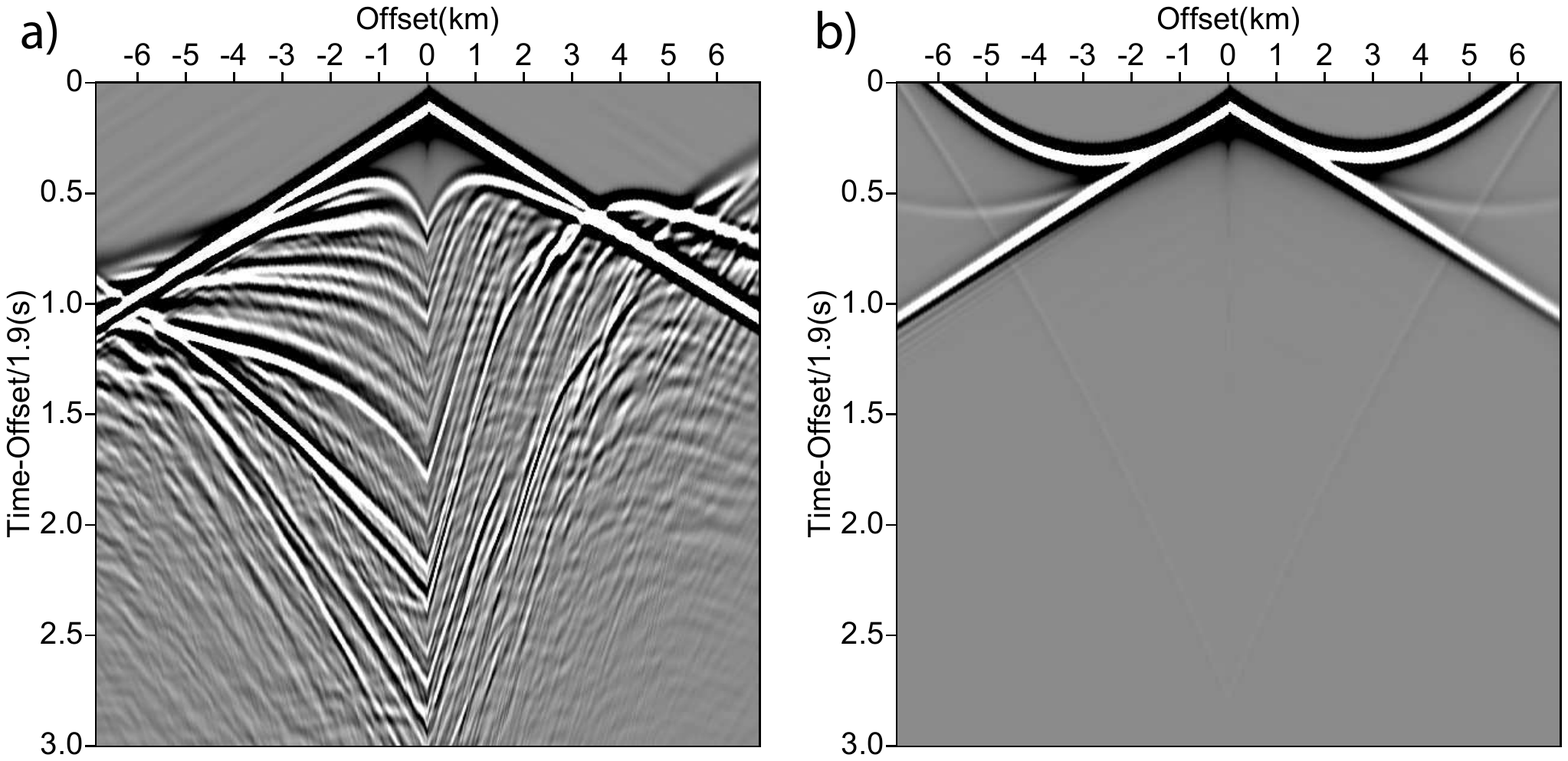}
\caption{Marmousi II example. Synthetic common-shot gather computed in the true velocity model (a) and the initial velocity model (b). The seismograms are plotted with a reduction velocity of
1900 m/s (origin time of each seismograms is $|\text{offset}|/1900$~s) to compress the time scale.}
\label{fig:fig_marmousi_true_init}
\end{center}
\end{figure}
A PML absorbing boundary condition is implemented on top of the grid (i.e., free-surface multiples are not involved in the inversion) and the source signature is a Ricker wavelet with a 10~Hz dominant frequency.
We design the inversion with a classical multiscale frequency continuation proceeding from the low frequencies to the higher ones.  One frequency is processed at a time between 3~Hz and 15~Hz with a 0.5~Hz interval, leading to 25 successive mono-frequency inversions. This frequency continuation strategy does not guarantee that, when the subsurface model is updated at a given frequency, the lower-frequency data are still fitted. To overcome this issue, we perform three cycles of multiscale inversion using the final model of the previous cycle as the initial model of the current cycle and starting the second and third cycles at 5~Hz and 7.5~Hz, respectively.  For each frequency, the stopping criterion of iteration is 
\begin{equation} \label{Stop}
k_{max}=10 \hspace{1 cm} \text{or}  \hspace{1 cm} (|| \bold{A(m}^k)\bold{u}^k-\bold{b}||_{F} \leq \delta  \hspace{1 cm} \text{and}  \hspace{1 cm} || \bold{Pu}^k-\bold{d}||_{F} \leq \epsilon_n), 
\end{equation}
with $\delta$=1e-3  and $\epsilon_n$=1e-5 and $k_{max}$ is the maximum number of iterations. Here, $F$ refers to the Frobenius norm. \\
The initial velocity model is a crude laterally-homogeneous velocity model in which the velocity linearly increases with depth (Figure \ref{fig:fig_marmousi_full_label}b). 
Comparison between a common-shot gather computed in the true and initial models highlights the kinematic inaccuracy of the initial velocity model, that would drive a classical FWI based
on a reduced-space formulation toward unacceptable local minimum (Figure~\ref{fig:fig_marmousi_true_init}). \\
A first set of inversions are performed with noiseless data. We use the same stopping criterion of iteration, equation~\ref{Stop}, and a constant penalty parameter $\lambda$ during iterations for WRI and IR-WRI. We set $\lambda$ equal to 1e-2 $\mu_1$. 
The final WRI and IR-WRI models are obtained after 615 and 533 iterations, respectively, and are shown in Figure~\ref{fig:fig_marmousi_full_label}(c-d). A direct comparison between the true model, the initial model and the two final WRI models along three vertical logs at horizontal distances of 4.5~km, 7~km and 9.5~km are shown in Figure~\ref{fig:fig_marmousi_log_label}a. These results show that IR-WRI successfully converges toward a satisfying model although a quite crude initial model. WRI leads to an acceptable velocity model in the shallow part down to roughly 2.5~km depth, with however a much poorer resolution level and less accurate positioning of sharp reflectors at depth than the model obtained with IR-WRI.  \\
%
We compute synthetic seismograms in the final WRI and IR-WRI models to assess how the errors in the velocity reconstruction translate into phase and amplitude data fit. We first superimpose in transparency the seismograms computed in the WRI and IR-WRI models on those computed in the true model to assess the traveltime match or in other words whether WRI and IR-WRI were affected by cycle skipping  (Figure~\ref{fig:fig_marmousi_direct}(a-b)).
This direct comparison suggests that none of the inversion was significantly hampered by cycle skipping, although we show small traveltime mismatches for deep reflections in the case of WRI (Figure~\ref{fig:fig_marmousi_direct}(a-b), black arrow). These traveltime mismatches can be related to the poor reconstruction of the deep reflectors above mentioned. Second, we show the residuals between the seismograms computed in the true model and the final WRI/IR-WRI models (Figure~\ref{fig:fig_marmousi_residuals}(a-b)).  These residuals clearly show that the iterative data and source residual updating allows for a better amplitude fit of both diving waves and reflections. This improved amplitude fit reflects the improved reconstruction of small-scale structures and velocity contrasts in the subsurface models when this updating is performed (Figures~\ref{fig:fig_marmousi_full_label}(c-d) and \ref{fig:fig_marmousi_log_label}a). Here, the improved resolution achieved by IR-WRI results because the iterative updating of the constraint violations allows us to 
fasten the convergence and, hence fulfill more closely the wave-equation and data-fitting constraints before the preset maximum number of iterations is reached, equation~\ref{Stop}. \\

%
\begin{figure}
\begin{center}
\includegraphics[scale=0.7]{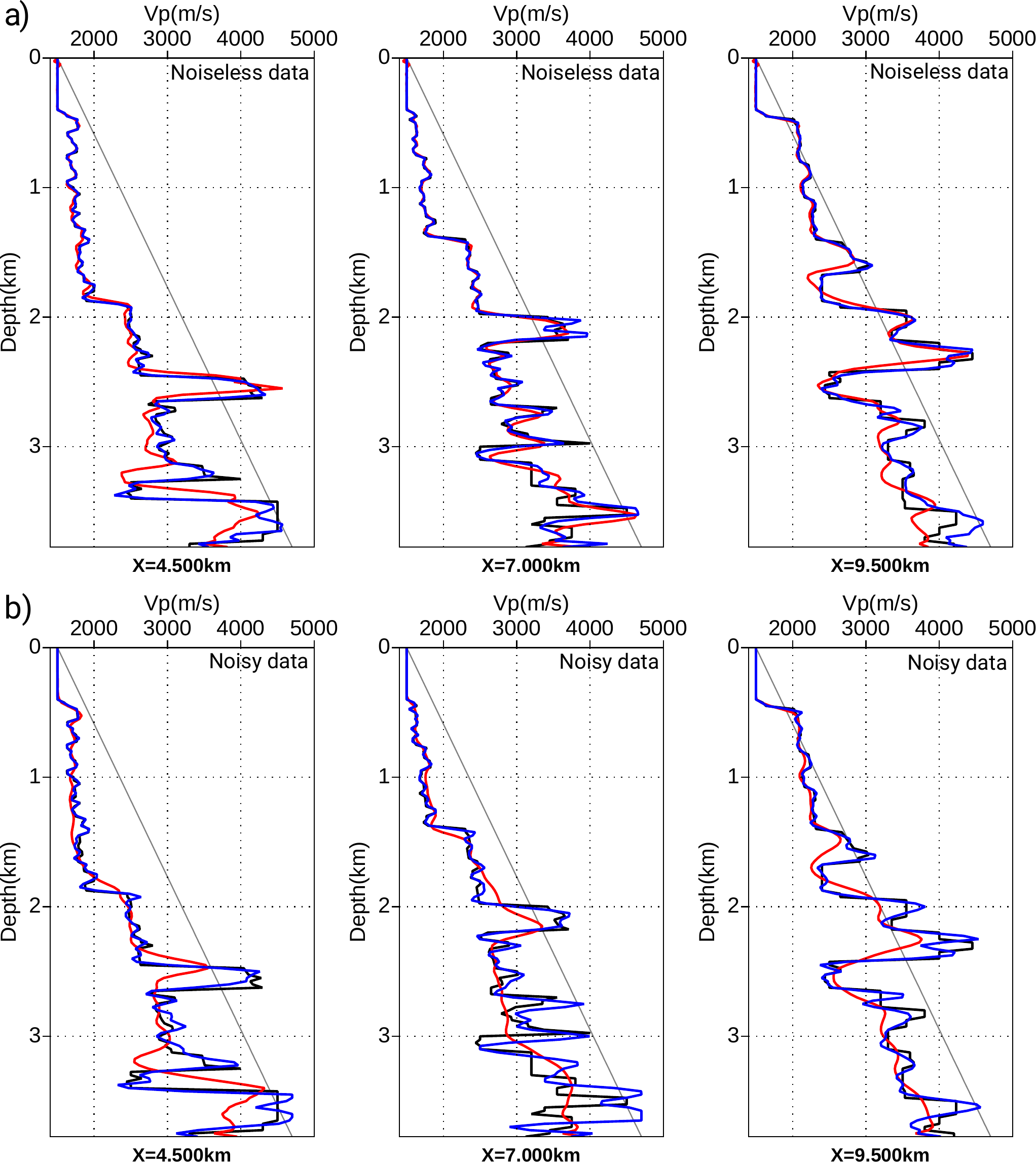}
\caption{Marmousi II example. Direct comparison between true (black), initial (gray), WRI (red) and IR-WRI (blue) velocity models along three logs at X=4.5km, 7km, 9.5km from left to right. (a) Noiseless data case (Figure~\ref{fig:fig_marmousi_full_label}c,d). (b) Noisy data case (Figure~\ref{fig:fig_marmousi_full_label}e,f). Note the deficit of resolution and mispositioning of deep reflectors in WRI results. These defaults increase significantly with noise, while IR-WRI reconstructions are more resilient to noise.}
\label{fig:fig_marmousi_log_label}
\end{center}
\end{figure}
To reproduce more realistic conditions, we add random noise to the recorded data with a signal-to-noise ratio (SNR) of 10~db. We added the noise in the frequency domain which means that each frequency involved in the inversion has the same SNR. We used the same experimental setup as above performing three inversion paths through the frequency batches. We used $\delta$=1e-3  and we set $\epsilon_n$ equal to the noise level of each frequency. The final WRI and IR-WRI models are obtained after 461 and 396 iterations, respectively, and are shown in Figure~\ref{fig:fig_marmousi_full_label}(e-f). A direct comparison between the true model, the initial model and the final WRI/IR-WRI models along the vertical logs are shown in Figure~\ref{fig:fig_marmousi_log_label}b. While the noise weekly impacts upon the IR-WRI results (compare Figures~\ref{fig:fig_marmousi_full_label}d and \ref{fig:fig_marmousi_full_label}f and Figures~\ref{fig:fig_marmousi_log_label}a, blue curves and \ref{fig:fig_marmousi_log_label}b, blue curves), the resolution of the model inferred from WRI is significantly degraded relative to that inferred from noiseless data, leading to a quite smooth reconstruction  (compare Figure~\ref{fig:fig_marmousi_full_label}c and \ref{fig:fig_marmousi_full_label}e and Figures~\ref{fig:fig_marmousi_log_label}a, red curves and \ref{fig:fig_marmousi_log_label}b, red curves). This highlights the higher resilience to noise of IR-WRI.

The data fit assessment for noisy data is shown in Figures~\ref{fig:fig_marmousi_direct}(c-d) and ~\ref{fig:fig_marmousi_residuals}(c-d)). The direct comparison between the seismograms computed in the true and WRI/IR-WRI models still reveal a good traveltime match without obvious evidence of cycle skipping (Figures~\ref{fig:fig_marmousi_direct}(c-d)). However, as for the noiseless case, we notice small traveltime mismatches in the case of WRI due to the poor focusing and mispositioning of deep reflectors shown in Figure~\ref{fig:fig_marmousi_log_label}a, red curves. More significantly, this direct comparison  shows that the WRI fails to match the amplitudes of several deep short-spread reflections due to poorly-resolved model reconstruction (Figures~\ref{fig:fig_marmousi_direct}c, white arrows), unlike IR-WRI (Figure~\ref{fig:fig_marmousi_direct}d). This is further confirmed by the data residuals shown in Figure~\ref{fig:fig_marmousi_residuals}(c-d). For both inversion, the amplitudes of the residuals as well as their low frequency content (due to the increase smoothness of the velocity reconstruction) have increased compared to those inferred from noiseless data. However, this increase is much more moderate in the case of IR-WRI.
%
%
\begin{figure}
\begin{center}
\includegraphics[scale=0.7]{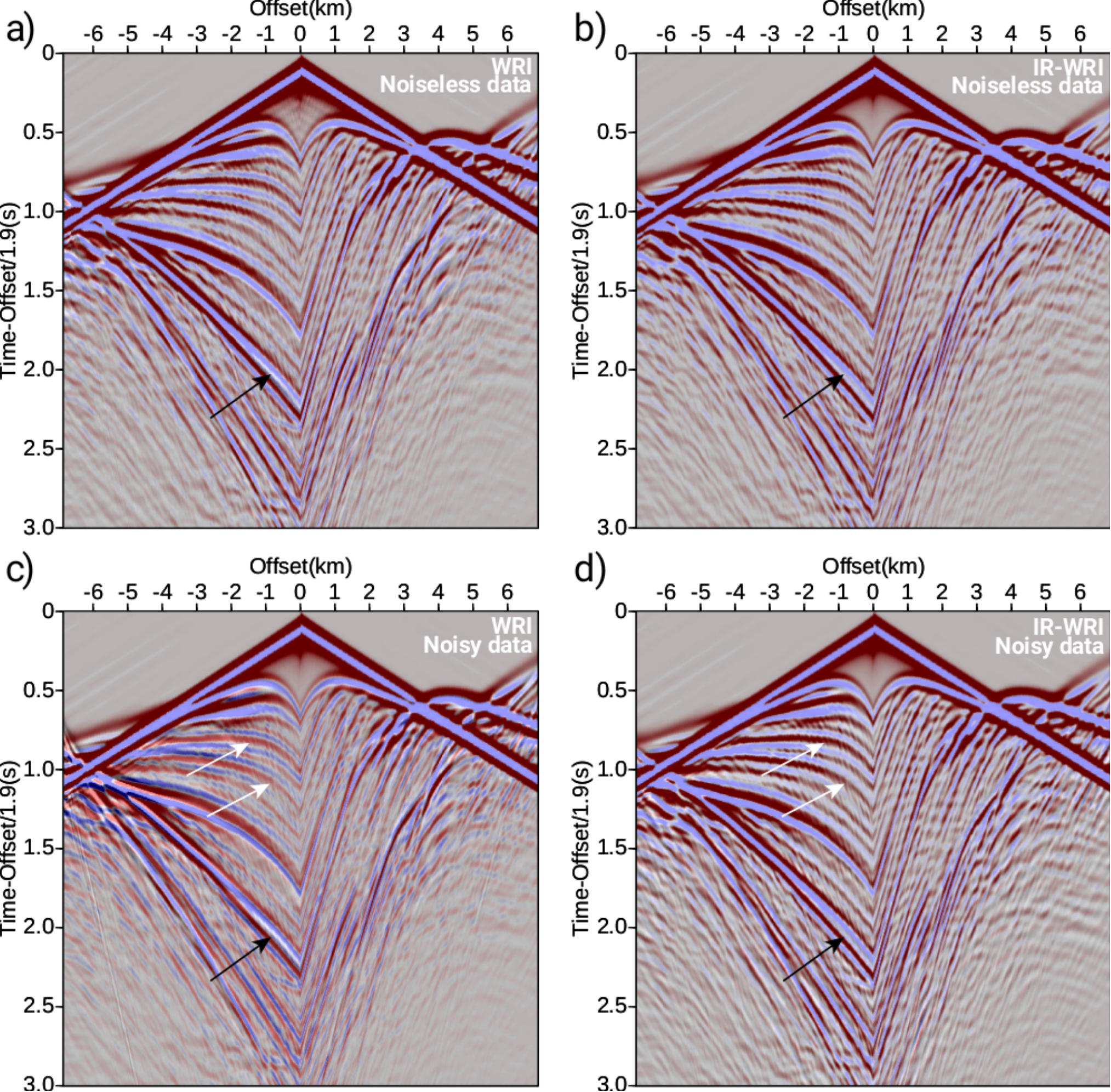}
\caption{Marmousi II example. Direct comparison between synthetic seismograms computed in true (red/white/blue scale scale) and WRI/IR-WRI (black/gray/white scale) models. Seismograms computed in the WRI/IR-WRI models are superimposed in transparency on seismograms computed in the true model. The two sets of seismograms are in phase if black/white wiggles of seismograms computed in the WRI/IR-WRI models match red/blue wiggles of seismograms computed in the true model. Seismograms are plotted with a reduction velocity of 1900 m/s. (a-b) For noiseless data, black/white seismograms are computed in (a) WRI and (b) IR-WRI models shown in Figure \ref{fig:fig_marmousi_full_label}c,d. (c-d) Same as (a-b) for noisy data (WRI/IR-WRI models are shown in Figures~\ref{fig:fig_marmousi_full_label}e,f). Arrows point the most obvious differences in the data fit achieved by WRI and IR-WRI. In the WRI results, data misfit take the form of significant amplitude underestimation (white arrows) and small time mismatches resulting from poorly resolved model reconstruction (black arrows).}
\label{fig:fig_marmousi_direct}
\end{center}
\end{figure}
%
%
\begin{figure}
\begin{center}
\includegraphics[scale=0.7]{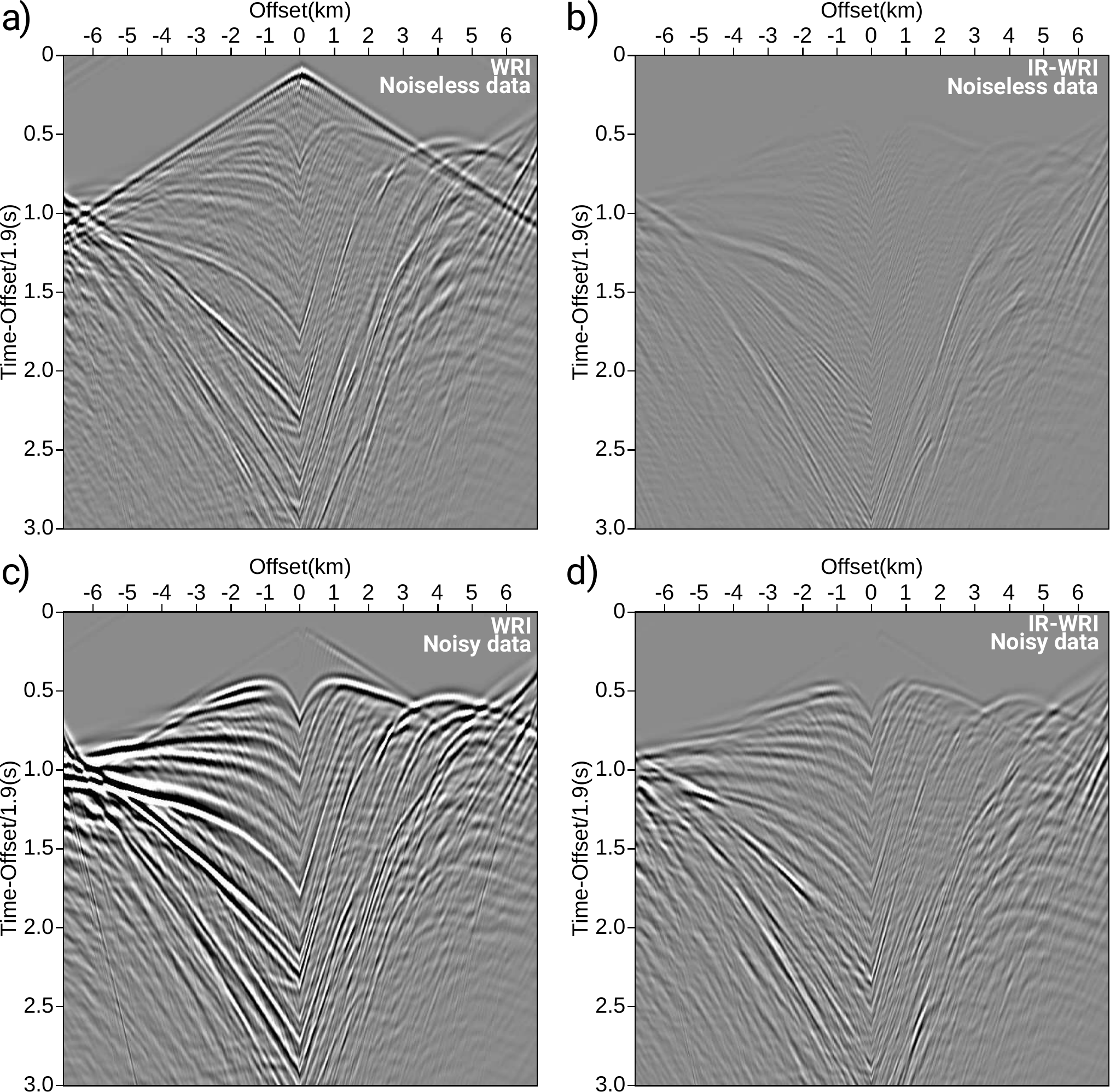}
\caption{Marmousi II example. Residuals between seismograms computed in the true and WRI/IR-WRI models. Each panel is plotted with the amplitude scale used in Figure~\ref{fig:fig_marmousi_true_init}.
(a-b) For noiseless data, residuals are computed from (a) WRI and (b) IR-WRI models. See Figure~\ref{fig:fig_marmousi_direct}(a-b) for a direct
comparison of the two sets of seismograms. (c-d) Same as (a-b) for noisy data. See Figure~\ref{fig:fig_marmousi_direct}(c-d) for a direct
comparison of the two sets of seismograms.}
\label{fig:fig_marmousi_residuals}
\end{center}
\end{figure}
%
\subsection{2004 BP salt model}

\begin{figure}
\begin{center}
\includegraphics[scale=0.7,trim=0cm 5cm 0cm 0cm]{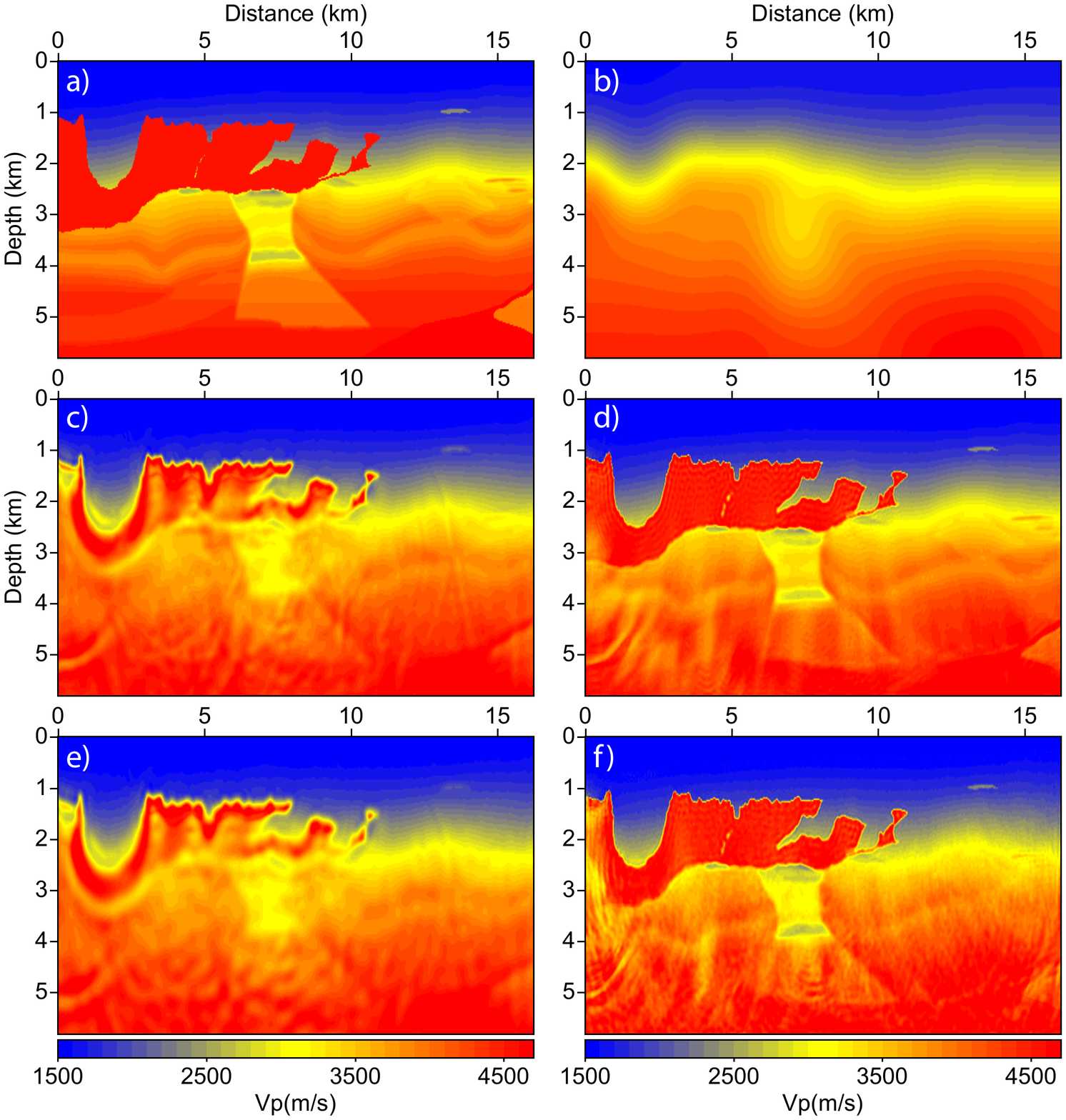}
\caption{2004 BP salt case study. (a) True 2004 BP model. (b) Initial velocity model. (c-d) Final WRI (c) and IR-WRI (d) velocity models for noiseless data. (e-f) Same as (c-d) for noisy data for a SNR of 10~db.}
\label{fig:fig_bp_full_label}
\end{center}
\end{figure}

Finally, we assess WRI and IR-WRI against a target of the challenging 2004 BP salt model (Figure \ref{fig:fig_bp_full_label}a) when a crude stating model is used (Figure \ref{fig:fig_bp_full_label}b). The selected target corresponds to the left part of the 2004 BP salt model and was  previously used in \citep{Metivier_2016_TOF} for an application of FWI based upon an optimal-transport distance. It is representative of the geology of the deep offshore Gulf of Mexico and mainly consists of a simple background with a complex rugose multi-valued salt body, sub-salt slow velocity anomalies related to over-pressure zones and a fast velocity anomaly to the right of the salt body \citep{Billette_2004_BPB}.   As for the Marmousi II model, we perform WRI and IR-WRI for noiseless and noisy data. The subsurface model is 16250~m wide and 5825~m deep and is discretized with a 25~m grid interval. We used 162 sources spaced 100~m apart on the top side of the model. The source signature is a Ricker wavelet with a 10~Hz dominant frequency. A line of receivers with a 25~m spacing are deployed at the surface leading to a stationary-receiver acquisition. As in the previous example, we used a 9-point anti-lumped mass staggered-grid stencil with PML boundary conditions for discretizing the Helmholtz equation. 
%
%
%
\begin{figure}
\begin{center}
\includegraphics[scale=0.7,trim=0cm 17cm 0cm 0cm]{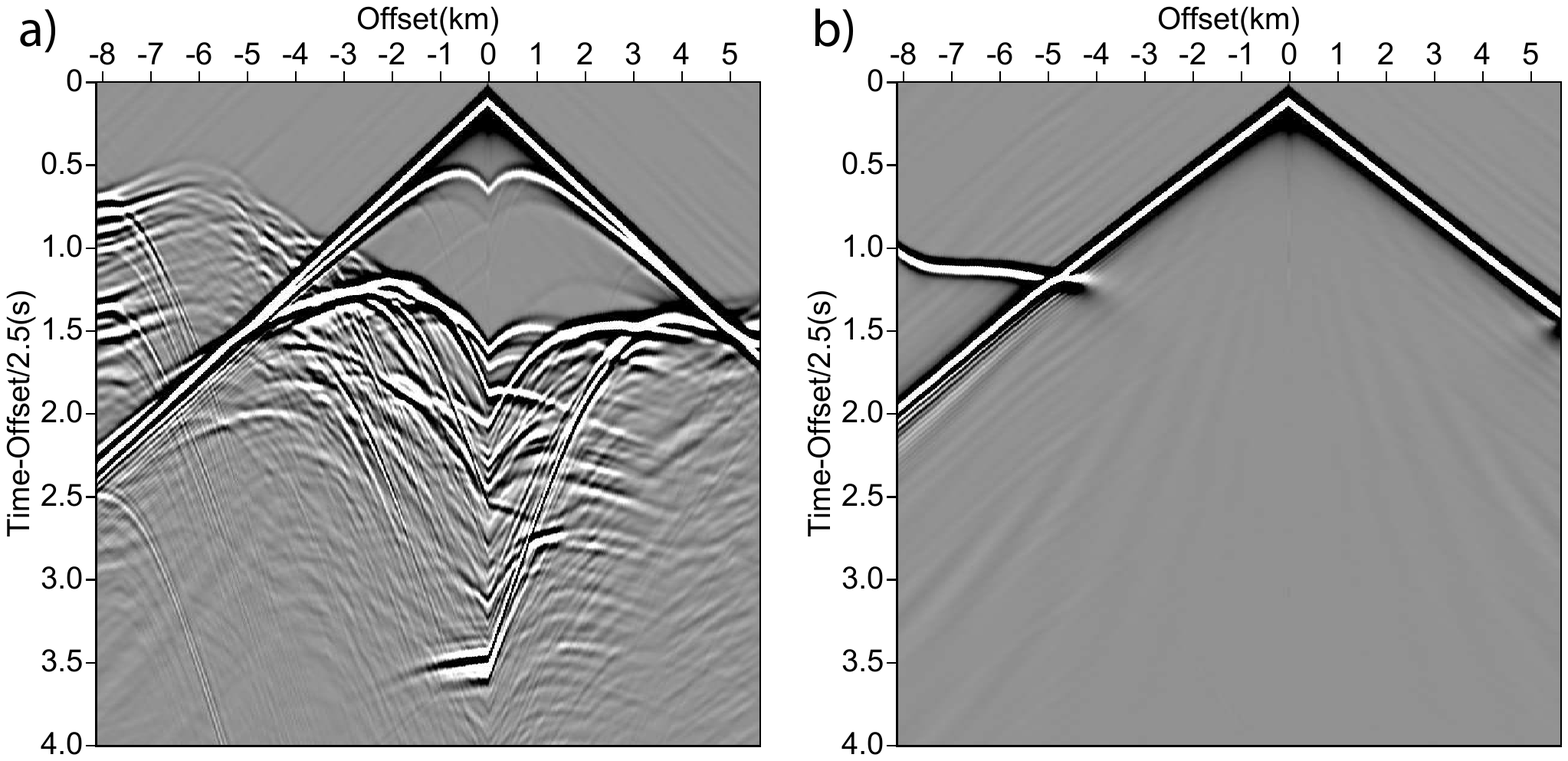}
\caption{2004 BP salt case study. Synthetic common-shot gather computed in the (a) true and (b) initial velocity models. Seismograms are plotted with a reduction velocity of
2500 m/s (the origin time of each seismograms is $|\text{offset}|/2500$~s) to compress the time scale.}
\label{fig:fig_bp_true_init}
\end{center}
\end{figure}
We used a smoothed version of the true velocity model as initial model (Figure \ref{fig:fig_bp_full_label}b). Following \citep{Metivier_2016_TOF}, we design the correlation length of the smoother such that the imprint of the salt body is canceled out. Comparison between a common-shot gather computed in the true and initial models shows that the footprint of the salt body in the wavefield has been removed and the diving waves are severely cycle skipped for a starting frequency of 3~Hz (Figure~\ref{fig:fig_bp_true_init}). 

\begin{figure}
\begin{center}
\includegraphics[scale=0.7]{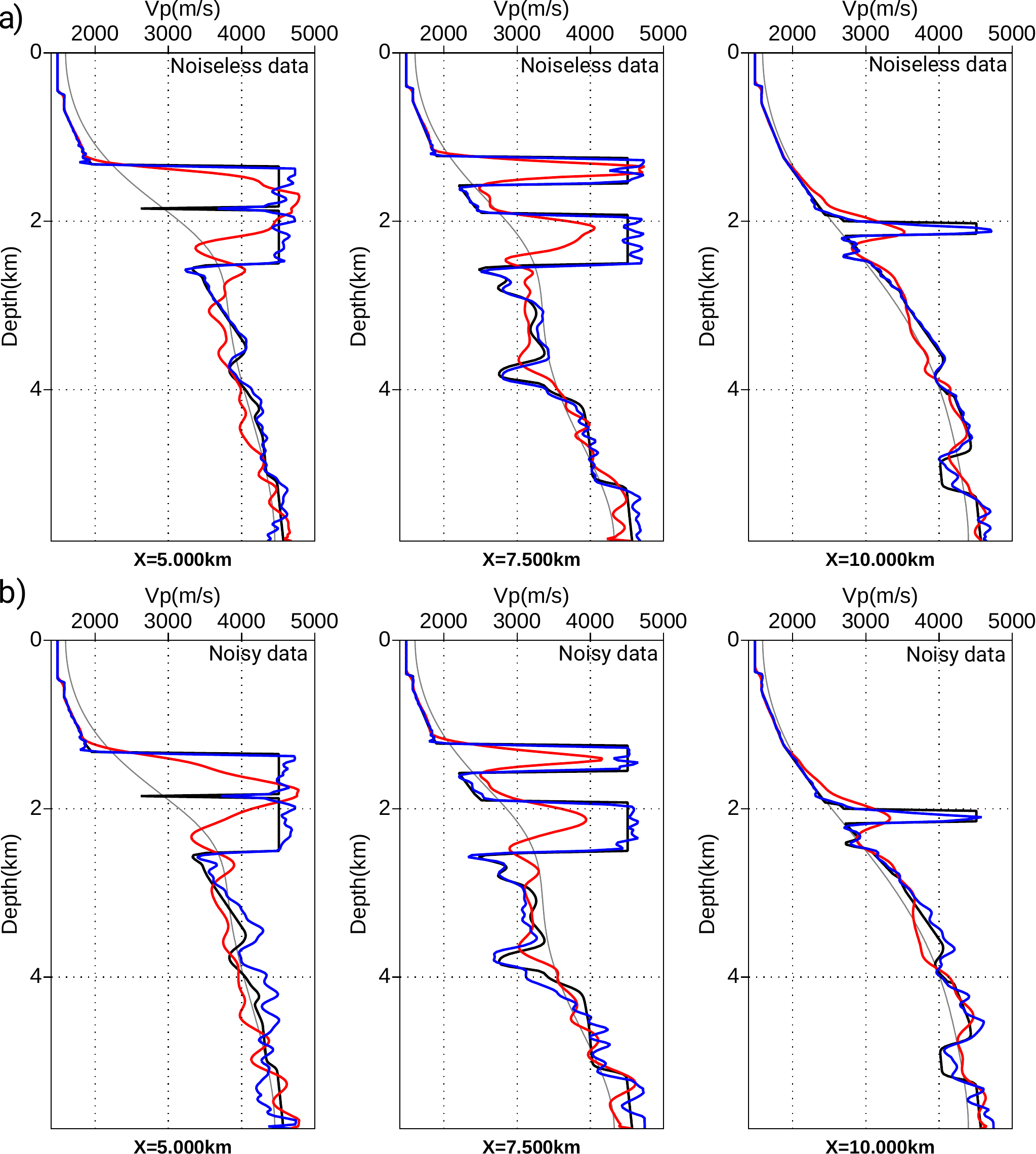}
\caption{2004 BP salt case study. Direct comparison between true (black), initial (gray), WRI (red) and IR-WRI (blue) velocity models along three logs at X=4.5km, 7km, 10km from left to right. (a) Noiseless data. (b) Noisy data.}
\label{fig:fig_bp_log_label}
\end{center}
\end{figure}

We used small batches of three frequencies with one frequency overlap between two consecutive batches, moving from the low frequencies to the higher ones according to a classical frequency continuation strategy. The starting and final frequencies are 3~Hz and 15~Hz and the sampling interval in one batch is 0.5~Hz. For each batch, the stopping criterion of iteration is given by equation \ref{Stop} with $k_{max}$=20. We use $\delta$=1e-3 for noiseless and noisy data, while $\epsilon_n$ is set to 1e-5 and the noise level of the batch for noiseless and noisy data, respectively. For the penalty parameter, we use $\lambda$=1e-3 $\mu_1$ and $\lambda$=2e-2 $\mu_1$ for noiseless and noisy data, respectively.
We perform three paths through the frequency batches to improve the inversion results (the starting frequency of the second and third path is 5~Hz and 7.5~Hz, respectively). 
The WRI and IR-WRI models inferred from noiseless data are shown in Figure~\ref{fig:fig_bp_full_label}(c-d). WRI and IR-WRI perform 561 and 448 iterations, respectively. Direct comparison between the true model, the starting model and the two WRI models along three vertical logs cross-cutting the salt body at 5~km, 7.5~km and 10~km distance is shown in Figure~\ref{fig:fig_bp_log_label}a.
We show that IR-WRI recovers with a quite high resolution the precise geometry of the rugose salt body, the sub-salt low-velocity anomalies as well as the shallow fast anomaly to the right (Figure~\ref{fig:fig_bp_full_label}d). In contrast, WRI successfully recovers the upper part of the salt but fails to retrieve the bottom part of the salt body. Moreover, the sub-salt structures are not recovered and the image of the fast anomaly to the right is blurred (Figure~\ref{fig:fig_bp_full_label}e). As for the Marmousi test, this results because the iterative updating of the data and wave-equation residuals contribute to fulfill more accurately the data-fitting and wave-equation constraints for $k_{max}$=20 due to faster convergence.

\begin{figure}
\begin{center}
\includegraphics[scale=0.7]{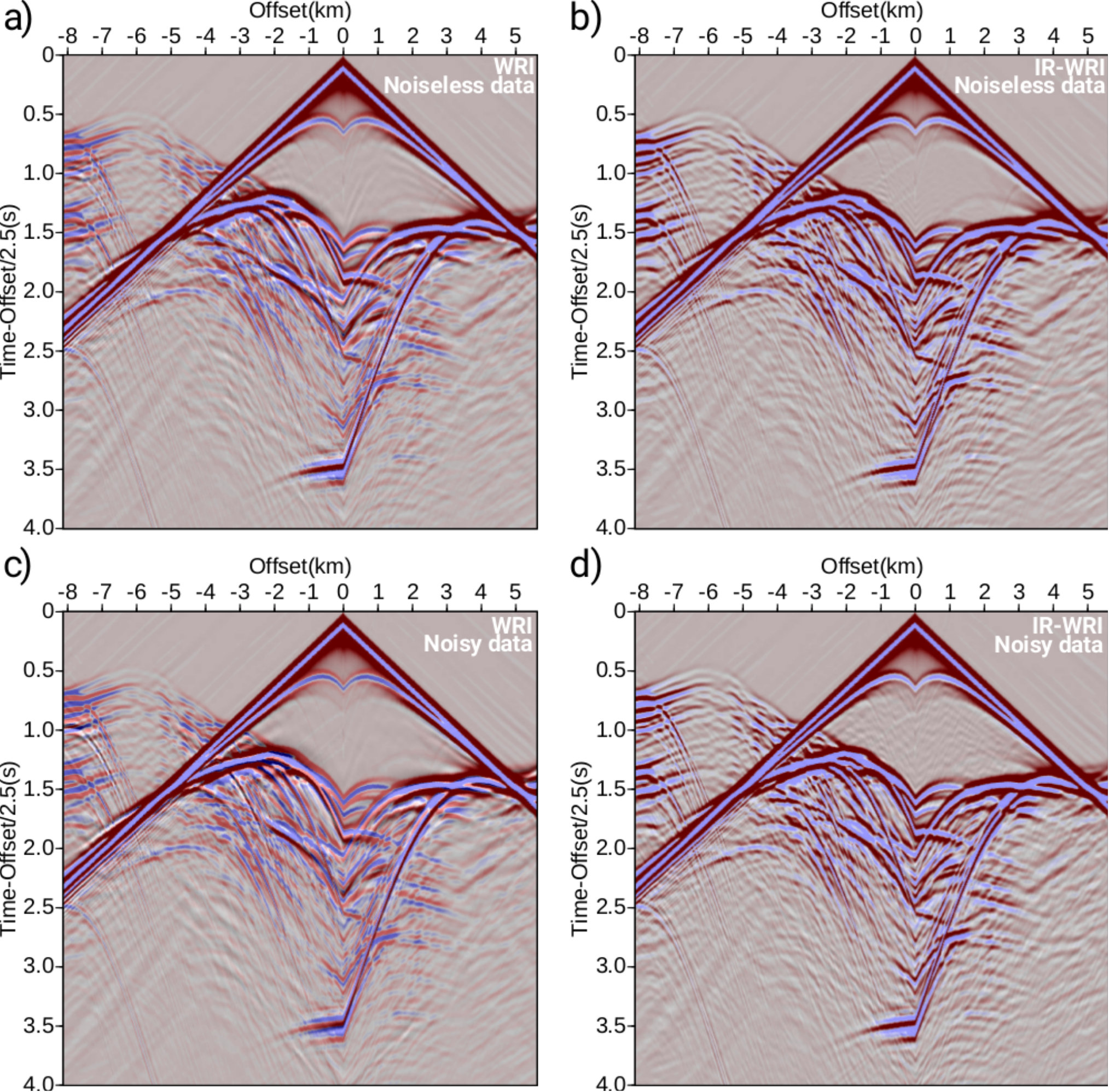}
\caption{2004 BP salt case study. Direct comparison between synthetic seismograms computed in true model (red/white/blue scale) and WRI/IR-WRI models (black/gray/white scale). Same showing as in Figure~\ref{fig:fig_marmousi_direct} is used to compare the two sets of seismograms. (a-b) For noiseless data, black/white seismograms are computed in (a) WRI and (b) IR-WRI models shown in Figure~\ref{fig:fig_bp_full_label}c,d. (c-d) Same as (a-b) for noisy data (WRI and IR-WRI models are shown in Figure~\ref{fig:fig_bp_full_label}e,f).  The arrows point the most obvious differences in the data fit achieved by WRI and IR-WRI. For WRI, data misfit take the form of significant amplitude underestimation and small time mismatch resulting from poorly resolved model reconstruction. Seismograms are plotted with a reduction velocity of 2500 m/s as in Figure~\ref{fig:fig_bp_true_init}.}
\label{fig:fig_bp_direct}
\end{center}
\end{figure}

\begin{figure}
\begin{center}
\includegraphics[scale=0.7]{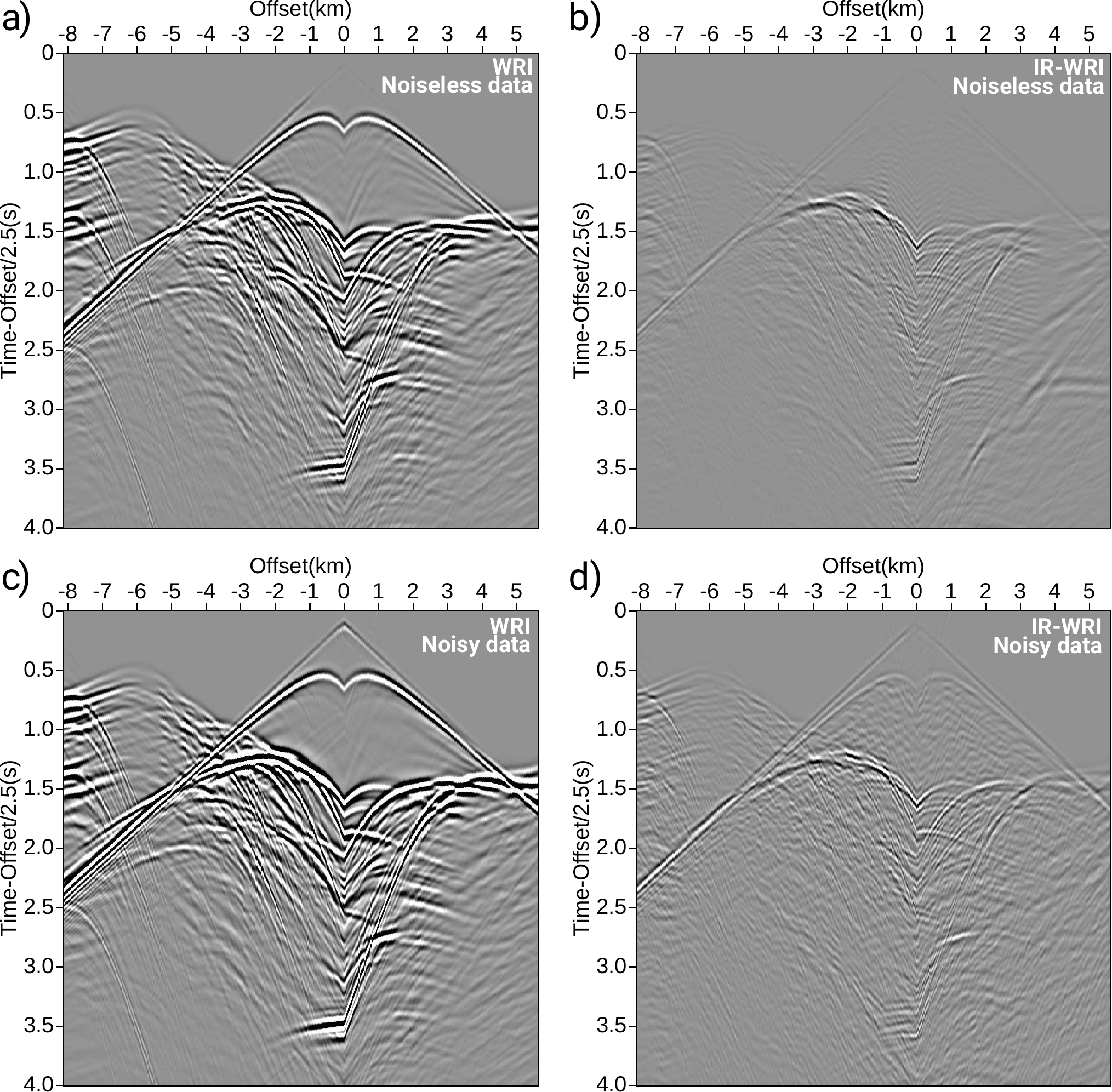}
\caption{2004 BP salt case study. Residuals between seismograms computed in true model and WRI/IR-WRI models. Each panel is plotted with the amplitude scale used in Figure~\ref{fig:fig_bp_true_init}.
(a-b) For noiseless data, residuals are computed from (a) WRI and (b) IR-WRI models. See Figure~\ref{fig:fig_bp_direct}(a-b) for a direct
comparison of the two sets of seismograms. (c-d) Same as (a-b) for noisy data. See Figure~\ref{fig:fig_bp_direct}(c-d) for a direct
comparison of the two sets of seismograms.}
\label{fig:fig_bp_residuals}
\end{center}
\end{figure}

As for the Marmousi example, the direct comparison between seismograms computed in the true and WRI/IR-WRI models does not show any evidence of cycle skipping (Figure~\ref{fig:fig_bp_direct}(a-b)). However, the residuals have much lower amplitudes in case of IR-WRI (Figure~\ref{fig:fig_bp_residuals}(a-b)) due to the improve reconstruction of the velocity contrasts and high wavenumbers highlighted in Figures~\ref{fig:fig_bp_full_label}d and \ref{fig:fig_bp_log_label}a.

When noisy data are used (Figures~\ref{fig:fig_bp_full_label}(e-f) and \ref{fig:fig_bp_log_label}b), the shallow imaging performed by the IR-WRI is weakly impacted upon by the noise. WRI and IR-WRI perform 461 and 396 iterations, respectively, and confirms again that the iterative residual updating not only improves the quality of the minimizer but fastens also the convergence. The sub-salt imaging is more affected by noise. However, the slow velocity anomalies below the salt body is still fairly well identified. In the WRI case, we show the same trend as the one revealed by the Marmousi example. The noise further degrades the results by hampering the reconstruction of short-scale features below the top of the salt but also in the overburden as revealed by the poor reconstruction of the fast anomaly on the right. 

The comparison of the seismograms computed in the true model and the WRI/IR-WRI models and the data residuals are shown in Figures~\ref{fig:fig_bp_direct}(c-d) and ~\ref{fig:fig_bp_residuals}(c-d).
A similar trend than for the Marmousi model is shown: When noise is added to the data, the amplitudes of the data residuals increase in much larger proportion in the WRI case.

%
%
%
%
%
%
%
\section{Discussion}
We have presented a significant yet easy-to-implement improvement of the WRI method of \citep{VanLeeuwen_2013_MLM}. 
We have first recast the original FWI problem as a feasibility problem where both data and source misfit are processed as constraints.
This feasibility problem is equivalent to constrained optimization problem where the objective function is identically zero. We solve this constrained problem with an augmented Lagrangian method (method of multiplier). A scaled form of the augmented Lagrangian shows that it is equivalent to a quadratic penalty function formed by the squares of the constraint violations (data and source misfit), to which are added the running sum of the data and source residuals of previous iterations. In the objectives of the penalty function, the data and the source represent the right-hand side of the constraints, and the running sum of the data and source residuals correspond to the scaled Lagrange multipliers. 
Following \citep{VanLeeuwen_2013_MLM}, we perform the refinement of the two primal variables (the wavefield and the subsurface parameters) in alternating mode, using the solution of one primal subproblem as a passive variable for the next subproblem. Compared to a classical application of ADMM where separable convex subproblems are solved in parallel, this sequential resolution of the two primal subproblems  allows us, on the one hand to linearize the parameter-estimation problem around the reconstructed wavefield and, on the other hand to manage the non-separability of the wavefield reconstruction and parameter estimation. This sequential solving of the two primal subproblems also prompts us to update the dual variable associated with the source residuals two times, once after each primal subproblem.

The convergence properties and the error bounds on the solution of penalty and augmented Lagrangian methods are discussed in the framework of nonlinear programming by \citep[][ Chapter 17]{Nocedal_2006_NOO}. As mentioned in the above reference (page 519), the augmented Lagrangian method gives us two ways of improving the accuracy of the minimizer: improve the accuracy of the Lagrange multiplier (namely, decrease the data and source residuals enough) or increase the penalty parameter, whereas the quadratic penalty approach gives us only one option: increase the penalty parameter. Accordingly and considering that we perform all of the numerical experiments with a fixed penalty parameter, the augmented-Lagrangian approach implemented in IR-WRI has logically converged to more accurate minimizers with a smaller number of iterations than the penalty method implemented in WRI.
The two leverages that control the accuracy of the IR-WRI minimizer also allow us to use moderate values of the penalty parameter without impacting prohibitively convergence speed. This moderate values allows the method to fit the data during the early iterations with large wave-equation error, which provides the most suitable framework to enlarge the search space and account for large time shifts, while satisfying the wave-equation constraint with small error at the convergence point. 
According to the analogy between WRI and IR-WRI highlighted by the scaled-form Lagrangian, IR-WRI can be viewed as a self-adaptive penalty method, where a tedious and potentially unstable continuous increasing of the penalty parameter is replaced by a stable dual steepest-ascent updating of the data and source residuals. This iterative residual updating in a quadratic misfit function is a well known procedure to refine solution of a wide class of linear inverse problems with gradient or Gauss-Newton steps as reminded in the appendix. 

The second key ingredient in IR-WRI is operator splitting implemented with an alternating direction strategy. IR-WRI relies on a generalization of the alternating-direction method of multiplier (ADMM) to biconvex problems \citep[][ section 9.2]{Boyd_2011_DOS}. A biconvex problem is one in which two variables can be partitioned into sets over which the problem is convex when the other variable is fixed \citep{Shen_2017_DMP}. IR-WRI is bi-convex because the wave-equation constraint is bilinear. This property leads to the linearization of the subproblem for $\bold{m}$ around $\bold{u}$ in our alternating-direction algorithm. Accordingly, the standard form of the ADMM can be readily implemented in IR-WRI, without however guarantee that IR-WRI benefits from the convergence properties of classical ADMM for convex problems.  The interested reader is refered to \citep{Bras_2012_ADA} for additional discussion on the applicability of ADMM on biconvex problems in the framework of elliptic PDE constraint. 

As a variant of WRI based upon the above-mention alternating direction strategy, \citep{vanLeeuwen_2016_PMP} recasted WRI as a reduced variable projection penalty method which leaves the parameter estimation subproblem nonlinear. \citep{Aravkin_2017_EQP} analyze the convergence properties of the reduced penalty method when the full Hessian is taken into account and conclude (their Theorem 2.5) that the variable projection penalty method is insensitive to the penalty parameter. More theoretical and numerical investigations are necessary to assess whether it is worth paying the price to solve a more complex nonlinear subproblem for $\bold{m}$ in the variable projection method compared to the ADMM algorithm implemented in IR-WRI which decomposes a complex problem into two simple ones. 

Although we have shown that the accuracy of the minimizer found by IR-WRI is weakly sensitive to the penalty parameter for a wide range of moderate values, the value which maximizes the convergence speed may be found accordingly to mitigate the computational cost. It is also worth reminding that updating the right-hand sides of the objectives in the IR-WRI penalty function does not introduce significant computational overhead in one iteration of IR-WRI compared to the WRI counterpart.

The reader familiar with the Bregman optimization method for constrained convex problem may have noticed that IR-WRI borrows some ingredients from this approach. The original Bregman method \citep{Bregman_1967_PCS} has been recast by \citep{Yin_2008_BIA,Goldstein_2009_SBM} as a penalty method where the right-hand side of the constraint is iteratively updated with the errors of previous iterations. 
Analogy between Bregman method and augmented Lagrangian method is recognized in \citep{Yin_2008_BIA}, while \citep{Esser_2009_ALA} reviewed the connection between ADMM and the split Bregman method.
The Bregman optimization method as well as its split Bregman extension has been mainly developed to solve l1-regularization problems for image denoising applications. A preliminary application of the split Bregman method to interface total-variation regularization with IR-WRI  is presented in \citep{Aghamiry_2018_MIA} where the identically-zero objective function in equation~\ref{main1_1} is replaced by the l1 norm of the total variation of the subsurface model. A short discussion on the choice of the penalty parameter in the framework of the Bregman optimization method  is provided in \citep[][ section 2.2]{Goldstein_2009_SBM} and is consistent with the results obtained with IR-WRI.

The original motivation of WRI was to make waveform inversion more resilient to cycle skipping. We apply WRI and IR-WRI on two complex synthetic examples to assess their ability to account for large time shifts.
The results show that both methods remain immune to cycle skipping for these case studies. However, only IR-WRI shows accurate reconstruction of deep targets such as the base of salt bodies and subsalt structures due to improved convergence history. This improved convergence manifests in the seismograms by a better amplitude fit of weak phases such as deep short-spread reflections. The improved convergence also makes IR-WRI more resilient to noise. These conclusions need to be validated against more realistic applications. Accordingly, perspective work aims to extend IR-WRI to 3D configuration and multiparameter reconstruction and assess the method on real stationary-receiver sea-bottom case studies as those tackled by \citep{Operto_2015_ETF,Operto_2018_MFF} where sharp contrasts generated by gas and chalky reservoir provide a suitable environment to assess the ability of IR-WRI to image large-contrast media 
and account for significant amplitude attenuation effects in realistic setting.

\section{Conclusion}
We have presented a new formulation of frequency-domain FWI based on wavefield reconstruction, which converges faster and improves the ability of the method to reconstruct short-scale structures, while preserving the resilience of the original formulation to cycle skipping. The improvement relies on an augmented Lagrangian method, which makes the method far less sensitive to the penalty parameter, when the value of this later is kept fixed in iterations. Accordingly, moderate values of the penalty parameter can be used such that the search space is significantly enlarged during the early iterations by relaxing the wave equation constraint without preventing to honor it with good accuracy at the final iteration. 
As the original penalty method, the augmented Lagrangian method is implemented with a splitting strategy, during which the wavefield reconstruction and the parameter estimation are performed in an alternating way and the solution of one subproblem is passed as
a passive variable for the next problem. This splitting strategy linearizes the parameter estimation around the reconstructed wavefield and manages the non-separability of the wavefield reconstruction and parameter estimation subproblems.
A scaled form of the augmented Lagrangian also shows that the new method is equivalent to a self-adaptive penalty method, where the dual updating of the Lagrange multipliers is recast as the iterative updating of the data and source right-hand sides in the objectives of the penalty function. Therefore, the new method does not generate significant computational overhead relative to the original penalty method. Numerical examples suggest that the iterative updating of both the data and source residuals dynamically manage the two competing objectives of the penalty function leading to improved convergence history. Application on the BP salt model reveals the ability of the method to reconstruct large-contrast structures including salt bodies and sub-salt structures, starting from crude initial models and a 3-Hz starting frequency. Synthetic seismogram modeling confirms that the inversion didn't suffer from cycle skipping for this case study and manages to reasonably fit amplitudes of weak phases such as deep short-spread reflections. 
Perspectives deal with multiparameter reconstruction, extension to 3D and application to real data.

\section{Acknowledgments}
This study was partially funded by the SEISCOPE consortium (\textit{http://seiscope2.osug.fr}), sponsored by AKERBP, CGG, CHEVRON, EXXON-MOBIL, JGI, PETROBRAS, SCHLUMBERGER, SHELL, SINOPEC, STATOIL and TOTAL. This study was granted access to the HPC resources of the Froggy platform of the CIMENT infrastructure (\textit{https://ciment.ujf-grenoble.fr}), which is supported by the Rh\^one-Alpes region (GRANT CPER07\_13 CIRA), the OSUG@2020 labex (reference ANR10 LABX56) and the Equip@Meso project (reference ANR-10-EQPX-29-01) of the programme Investissements d'Avenir supervised by the Agence Nationale pour la Recherche, and the HPC resources of CINES/IDRIS/TGCC under the allocation 046091 made by GENCI."

\appendix
\section{\\Iterative solution refinement for linear inverse problems}
\label{appendixa}
When one solves a system of equations, the computed solution may deviate from the desired one. In linear algebra, this might result from round-off errors in large-scale problems \citep[e.g.][ page 61] {Press_2007_NR3}. In convex and nonlinear optimization, this might arise when regularization or penalty methods are used to solve an ill-conditioned problem \citep{Gholami_2016_RGI}. Another example deals with
approximate solution of nonlinear inverse problem through linearization, as for example least-squares migration/inversion, where the Born approximation, sometimes combined with ray theory, is used to linearize the wave equation around a background velocity model \citep{Lambare_1992_IAI,Jin_1992_TDA} or linearized AVO inversion \citep{Gholami_2017_CNA}. The residual errors in the computed solution may require to solve repeatedly the problem to iteratively refine the solution accuracy. In this appendix, we review this iterative solution refinement procedure which amounts to update at a given iteration the right-hand side of the linear problem with the running sum of the residuals of previous iterations. 

Let's solve iteratively
\begin{equation}
\bold{A} \bold{x} = \bold{b}.
\label{eqa1}
\end{equation}
At iteration 1, we find an approximate solution $\bold{x}_1$ by solving the quadratic program in a least-squares sense
\begin{equation}
\bold{x}_1 = \text{arg} \min_x \| \bold{A} \bold{x} - \bold{b} \| = \bold{A}^{-g} \bold{b}.
\label{eqinit}
\end{equation}
where $\bold{A}^{-g}$ denotes the generalized or pseudo inverse inverse of $\bold{A}$, $\bold{A}^{-g} = \left( \bold{A}^T \bold{A} + \beta \bold{I}\right)^{-1} \bold{A}^T$ and $\beta$ is a stabilizing parameter \citep{Menke_2012_GDA}. \\
The solution $\bold{b}_1$ of the forward problem is given by
\begin{equation}
\bold{b}_1 = \bold{A}  \bold{x}_1
\label{eqa3}
\end{equation}
and fits $\bold{b}$ with an error $\delta \bold{b}_1 = \bold{b} - \bold{b}_1$.
%
%
We search the perturbation $\delta \bold{x} = \bold{x} - \bold{x}_1$ that needs to be applied on $\bold{x}_1$ to correct the residual $\delta \bold{b}_1$:
\begin{equation}
\delta \bold{x} = \text{arg} \min_x \| \bold{A} \delta \bold{x} - \delta \bold{b}_1 \| = \bold{A}^{-g} \delta \bold{b}_1 = \bold{A}^{-g} \left(\bold{b} - \bold{A} \bold{x}_1\right).
\label{eqmin}
\end{equation}
The refined solution is given by
\begin{equation}
\bold{x}_2 = \bold{x}_1 + \bold{A}^{-g} \delta \bold{b}_1  = \bold{A}^{-g} \bold{b} + \bold{A}^{-g} \left(\bold{b} - \bold{A} \bold{x}_1\right).
\end{equation}
Proceeding with the next iteration leads to the following refined solution
\begin{equation}
\bold{x}_3 = \bold{x}_2 + \bold{A}^{-g} (\bold{b} - \bold{A} \bold{x}_2) =  \bold{A}^{-g} \bold{b} + \bold{A}^{-g} \left(\bold{b} - \bold{A} \bold{x}_1\right) + \bold{A}^{-g} (\bold{b} - \bold{A} \bold{x}_2).
\end{equation}
Repeating this procedure $k$ times, we obtain
\begin{equation}
\bold{x}_k =  \bold{A}^{-g} \bold{b} + \sum_{i=1}^{k-1} \bold{A}^{-g}  \left(\bold{b} - \bold{A} \bold{x}_i\right) = \bold{A}^{-g} \left( \bold{b} +  \sum_{i=1}^{k-1}  (\bold{b} - \bold{A} \bold{x}_i) \right).
\label{eqa10}
\end{equation}
Analogy between equations~\ref{eqinit} and \ref{eqa10} allows us to write $\bold{x}_k$ as the solution of the quadratic program
\begin{equation}
\bold{x}_k = \text{arg} \min_x \| \bold{A} \bold{x} - \bold{b} - \Delta \bold{b}_{k-1} \|,
\label{eqmin1}
\end{equation}
where the original right-hand side $\bold{b}$ has been updated with the running sum of the residuals of the first $k-1$ iterations, $\Delta \bold{b}_{k-1}=  \sum_{i=1}^{k-1}  (\bold{b} - \bold{A} \bold{x}_i)$.
This means that, instead of looking for a solution perturbation and adding this perturbation to the current solution to refine it, the imprint of the current solution is directly injected in the objective via its 
linear relationship with the residuals of previous iterations.
This iterative updating of the right-hand sides of a linear problem describes the approach that is used in this study to refine the solutions of the wavefield reconstruction and subsurface model
estimation subproblems when the initial Lagrangian method is recast as a penalty method after scaling the dual variables by the penalty parameter, equation~\ref{ALS}.

\bibliographystyle{seg}
\bibliography{../../BIBLIO/biblioseiscope,../../BIBLIO/bibliotmp}

\end{document}